\documentclass[12pt, a4paper]{article}

\usepackage{fixltx2e}
\usepackage{cmap}
\usepackage[utf8]{inputenc}
\usepackage[T1]{fontenc}
\usepackage[english]{babel}
\usepackage{amsmath, amsfonts, amssymb, amsthm, mathtools, indentfirst, bbm}
\usepackage[shortlabels]{enumitem}
\usepackage[languagenames,fixlanguage]{babelbib}
\usepackage{a4wide}
\allowdisplaybreaks
\usepackage{commath}
\usepackage{textcomp}

\usepackage{setspace}
\usepackage[unicode=true, pdfstartview={XYZ null null 1.00}, colorlinks=true, linkcolor=black, urlcolor=black, citecolor=black]{hyperref}

\newcommand{\re}{\operatorname{Re}}

\newcommand{\bbN}{\mathbb{N}}
\newcommand{\bbZ}{\mathbb{Z}}
\newcommand{\bbR}{\mathbb{R}}

\newcommand{\bbP}{\mathbb{P}}
\newcommand{\bbE}{\mathbb{E}}
\newcommand{\bbjedan}{\mathbbm{1}}

\newcommand{\calL}{\mathcal{L}}
\newcommand{\calF}{\mathcal{F}}

\newcommand{\aps}[1]{\vert #1 \vert}

\newcommand{\OBL}[1]{\left( #1 \right)}
\newcommand{\UGL}[1]{\left[ #1 \right]}
\newcommand{\VIT}[1]{\left\{ #1 \right\}}

\newcommand{\FJADEF}[3]{{#1}:{#2}\to{#3}}

\let\ge\geqslant
\let\le\leqslant

\newtheorem{DEF}{Definition}[section]
\newtheorem{REM}[DEF]{Remark}
\newtheorem{TM}[DEF]{Theorem}
\newtheorem{PROP}[DEF]{Proposition}
\newtheorem{LM}[DEF]{Lemma}
\newtheorem{COR}[DEF]{Corollary}

\makeatletter
\newenvironment{PF}[1][\proofname]{\par
	\pushQED{\qed}%
	\normalfont \topsep0\p@\relax
	\trivlist
	\item[\hskip\labelsep\bfseries
	#1\@addpunct{:}]\ignorespaces
}{%
	\popQED\endtrivlist\@endpefalse
	\addvspace{2\topsep}
}
\makeatother
\renewcommand{\qed}{\hfill $\blacksquare$}

\numberwithin{equation}{section}

\title{Harnack inequality for subordinate random walks}
\author{Ante Mimica ($\dagger$) \and Stjepan Šebek}
\date{}
\begin{document}
\maketitle

\begin{abstract}
\noindent
In this paper, we consider a large class of subordinate random walks $X$ on integer lattice $\mathbb{Z}^d$ via subordinators with Laplace exponents which are complete Bernstein functions satisfying a certain lower scaling condition at zero. We establish estimates for one-step transition probabilities, the Green function and the Green function of a ball, and prove the Harnack inequality for non-negative harmonic functions.

\vspace{15px}


\end{abstract}

\section{Introduction}
\noindent
Let $(Y_k)_{k \ge 1}$ be a sequence of independent, identically distributed random variables defined on a probability space $(\Omega, \calF, \bbP)$, taking values in the integer lattice $\bbZ^d$, with distribution $\bbP(Y_k = e_i) = \bbP(Y_k = -e_i) = 1/2d$, $i = 1, 2, \ldots, d$, where $e_i$ is the $i$-th vector of the standard basis for $\bbR^d$. A simple symmetric random walk in $\bbZ^d$ $(d \ge 1)$ starting at $x \in \bbZ^d$ is a stochastic process $Z = (Z_n)_{n \ge 0}$, with $Z_0 = x$ and $Z_n = x + Y_1 + \cdots + Y_n$.

Let $Z = (Z_n)_{n \ge 0}$ be a simple symmetric random walk in $\bbZ^d$ starting at the origin. Further, let
\begin{equation*}
	\phi(\lambda) := \int_{\langle 0, \infty \rangle}{\OBL{1 - e^{-\lambda t}} \mu (\dif t)}
\end{equation*}
be a Bernstein function satisfying $\phi(1)=1$. Here $\mu$ is a measure on $\langle 0,\infty \rangle$ satisfying $\int_{\langle 0,\infty \rangle} {(1 \wedge t) \mu(\dif t)} < \infty$ called the L\'{e}vy measure. For $m \in \bbN$ denote
\begin{equation}\label{eq:def_of_c_m^phi}
	c_m^{\phi} := \int_{\langle 0, \infty \rangle}{\frac{t^m}{m!}e^{-t}\mu (dt)}
\end{equation}
and notice that
\begin{equation*}
	\sum_{m = 1}^{\infty} {c_m^{\phi}} = \int_{\langle 0, \infty \rangle} {(e^t - 1) e^{-t} \mu(\dif t)} = \int_{\langle 0, \infty \rangle} {(1 - e^{-t}) \mu(\dif t)} = \phi(1) = 1.
\end{equation*}
Hence, we can define a random variable $R$ with $\bbP(R = m) = c_m^{\phi}$, $m \in \bbN$. Now we define the random walk $T = (T_n)_{n \ge 0}$ on $\bbZ_+$ by $T_n := \sum_{k = 1}^{n} {R_k}$, where $(R_k)_{k \ge 1}$ is a sequence of independent, identically distributed random variables with the same distribution as $R$ and independent of the process $Z$. Subordinate random walk is a stochastic process $X = (X_n)_{n \ge 0}$ which is defined by $X_n := Z_{T_n}$, $n \ge 0$. It is straightforward to see that the subordinate random walk is indeed a random walk. Hence, there exists a sequence of independent, identically distributed random variables $(\xi_k)_{k \ge 1}$ with the same distribution as $X_1$ such that
\begin{equation}\label{eq:X_as_RW}
	X_n \overset{d}{=} \sum_{k = 1}^n {\xi_k}, \quad n \ge 0.
\end{equation}
We can easily find the explicit expression for the distribution of the random variable $X_1$:
\begin{align}\label{eq:definition_of_xi}
	\bbP(X_1 = x)
	& =  \bbP(Z_{T_1} = x) = \bbP(Z_{R_1} = x) = \sum_{m = 1}^{\infty} {\bbP(Z_{R_1} = x \mid R_1 = m)} c_m^{\phi} \nonumber \\
	& =  \sum\limits_{m = 1}^{\infty} \int_{\langle 0, \infty \rangle}{\frac{t^m}{m!} e^{-t} \mu(dt)} \bbP(Z_m = x), \quad x \in \bbZ^d.
\end{align}
We denote the transition matrix of the subordinate random walk $X$ with $P$, i.e.\ $P = (p(x, y) : x, y \in \bbZ^d)$, where $p(x, y) = \bbP(x + X_1 = y)$.

We will impose some additional constraints on the Laplace exponent $\phi$. First, $\phi$ will be a complete Bernstein function \cite[Definition 6.1.]{SS12} and it will satisfy the following lower scaling condition: there exist $0 < \gamma_1 < 1$ and $a_1 > 0$ such that
\begin{equation}\label{eq:svojstvo_skaliranja}
	\frac{\phi(R)}{\phi(r)} \ge a_1 \OBL{\frac{R}{r}}^{\gamma_1}, \quad \forall\,\, 0 < r \le R\le 1.
\end{equation}
In dimension $d \le 2$, we additionally assume that there exist $\gamma_1 \le \gamma_2 < 1$ and $a_2 > 0$ such that
\begin{equation}\label{eq:upper_svojstvo_skaliranja}
	\frac{\phi(R)}{\phi(r)} \le a_2 \OBL{\frac{R}{r}}^{\gamma_2}, \quad \forall\,\, 0 < r \le R\le 1.
\end{equation}
It is well known that,  if $\phi$ is a Bernstein function, then $\phi(\lambda t) \le \lambda \phi(t)$ for all $\lambda \ge 1$, $t > 0$, implying $\phi(v)/v \le \phi(u)/u$, $0 < u \le v$. Using these two facts, proof of the next lemma is straightforward.
\begin{LM}\label{lm:free_upper_scaling}
	If $\phi$ is a Bernstein function, then for all $\lambda, t > 0$, $1 \wedge \lambda \le \phi(\lambda t) / \phi(t) \le 1 \vee \lambda$.
\end{LM}
Using Lemma \ref{lm:free_upper_scaling} we get $\phi(R) / \phi(r) \le R/r$ for all $0 < r \le R \le 1$ and this suffices for $d \ge 3$.

The main result of this paper is a scale invariant Harnack inequality for subordinate random walks. The proof will be given in the last section. Before we state the result, we define the notion of harmonic function with respect to subordinate random walk $X$.
\begin{DEF}\label{def:harmonic_function}
	We say that a function $\FJADEF{f}{\bbZ^d}{[0, \infty\rangle}$ is harmonic in $B \subset \bbZ^d$, with respect to $X$, if
	\begin{equation*}
		f(x) = Pf(x) = \sum_{y \in \bbZ^d} {p(x, y) f(y)}, \quad \forall\ x \in B.
	\end{equation*}
\end{DEF}
This definition is equivalent to the mean-value property in terms of the exit from a finite subset of $\bbZ^d$: If $B \subset \bbZ^d$ is finite then $\FJADEF{f}{\bbZ^d}{[0, \infty \rangle}$ is harmonic in $B$, with respect to $X$, if and only if $f(x) = \bbE_x[f(X_{\tau_B})]$ for every $x \in B$, where $\tau_B := \inf \{n \ge 1 : X_n \notin B\}$.

For $x \in \bbZ^d$ and $r > 0$ we define $B(x, r) := \{y \in \bbZ^d : \aps{y - x} < r\}$. We use shorthand notation $B_r$ for $B(0, r)$.
\begin{TM}[Harnack inequality]\label{tm:Harnack_inequality}
	Let $X = (X_n)_{n \ge 0}$ be a subordinate random walk in $\bbZ^d$, $d \ge 1$, with $\phi$ a complete Bernstein function satisfying \eqref{eq:svojstvo_skaliranja} and in the case $d \le 2$ also \eqref{eq:upper_svojstvo_skaliranja}. For each $a < 1$, there exists a constant $c_a < \infty$ such that if $\FJADEF{f}{\bbZ^d}{[0, \infty\rangle}$ is harmonic on $B(x, n)$, with respect to $X$, for $x \in \bbZ^d$ and $n \in \bbN$, then
	\begin{equation*}
		f(x_1) \le c_a f(x_2), \quad \forall\ x_1, x_2 \in B(x, an).
	\end{equation*}
\end{TM}
Notice that the constant $c_a$ is uniform for all $n \in \bbN$. That is why we call this result scale invariant Harnack inequality.

Some authors have already dealt with this problem and Harnack inequality was proved for symmetric simple random walk in $\bbZ^d$ \cite[Theorem 1.7.2]{La96} and random walks with steps of infinite range, but with some assumptions on the Green function and some restrictions such as finite second moment of the step \cite{BK08, LP93}.

Notion of discrete subordination was developed in \cite{BS12} and it was discussed in details in \cite{BC10}, but under different assumptions on $\phi$ than the ones we have. Using discrete subordination we can obtain random walks with steps of infinite second moment, see Remark \ref{rem:infinite_variance_of_the_step}. Harnack inequality has not been proved so far for such random walks.

In Section \ref{sec:Preliminaries} we state an important result about gamma function that we use later, we discuss under which conditions subordinate random walk is transient and we introduce functions $g$ and $j$ and examine their properties. The estimates of one-step transition probabilities of subordinate random walk are given in Section \ref{sec:Transition_probability_estimates}. In Section \ref{sec:Green_function_estimates} we derive estimates for the Green function. This is very valuable result which gives the answer to the question posed in \cite[Remark 1]{BCT15}. Using estimates developed in Section \ref{sec:Transition_probability_estimates} and \ref{sec:Green_function_estimates} and following \cite[Section 4]{MK12}, in Section \ref{sec:Estimates_of_Green_function_of_a_ball} we find estimates for the Green function of a ball. In Section \ref{sec:Proof_of_Harnack_inequality} we introduce Poisson kernel and prove Harnack inequality.

Throughout this paper, $d\ge 1$ and the constants $a_1, a_2$, $\gamma_1, \gamma_2$ and $C_i$, $i = 1, 2, \ldots, 9$ will be fixed. We use $c_1, c_2, \ldots$ to denote generic constants, whose exact values are not important and can change from one appearance to another. The labeling of the constants $c_1, c_2, \ldots$ starts anew in the statement of each result. The dependence of the constant $c$ on the dimension $d$ will not be mentioned explicitly. We will use ``:='' to denote a definition, which is read as ``is defined to be''. We will use $\dif x$ to denote the Lebesgue measure in $\bbR^d$. We denote the Euclidean distance between $x$ and $y$ in $\bbR^d$ by $\aps{x - y}$. For $a, b \in \bbR$, $a \wedge b := \min\{a, b\}$ and $a \vee b := \max\{a, b\}$. For any two positive functions $f$ and $g$, we use the notation $f \asymp g$, which is read as ``$f$ is comparable to $g$'', to denote that there exist some constants $c_1, c_2 > 0$ such that $c_1 f \le g \le c_2 f$ on their common domain of definition. We also use notation $f \sim g$ to denote that $\lim_{x \to \infty} f(x) / g(x) = 1$.

\section{Preliminaries}\label{sec:Preliminaries}
\noindent
In this section we first state an important result about the ratio of gamma functions that is needed later. Secondly, we discuss under which conditions subordinate random walk is transient. At the end of this section we define functions $g$ and $j$ that we use later and we prove some of their properties.

\subsection{Ratio of gamma functions}
\begin{LM}\label{lm:incomplete_complete_gamma_ratio}
	Let $\Gamma(x, a) = \int_a^{\infty} {t^{x - 1} e^{-t} \dif t}$. Then		\begin{equation}\label{eq:incomplete_complete_gamma_ratio}
		\lim_{x \to \infty} {\frac{\Gamma(x + 1, x)}{\Gamma(x + 1)}} = \frac{1}{2}.
	\end{equation}
\end{LM}
\begin{PF}
	Using a well-known Stirling's formula
		\begin{equation}\label{eq:stirlings_formula}
	\Gamma(x + 1) \sim \sqrt{2 \pi x}\ x^x e^{-x}, \quad x \to \infty
	\end{equation}
	and \cite[Formula 6.5.35]{AS72} that states
	\begin{equation*}
	 	\Gamma(x + 1, x) \sim \sqrt{2^{-1} \pi x} \ x^x e^{-x}, \quad x \to \infty
\end{equation*}
	we get
	\begin{equation*}
		\lim_{x \to \infty} {\frac{\Gamma(x + 1, x)}{\Gamma(x + 1)}} = \lim_{x \to \infty} \frac{\sqrt{2^{-1} \pi x}\ x^x e^{-x}}{\sqrt{2 \pi x} \ x^x e^{-x}} = \frac{1}{2}.
	\end{equation*}
\end{PF}

\subsection{Transience of subordinate random walks}
\noindent
Our considerations only make sense if the random walk that we defined is transient. In the case of a recurrent random walk, the Green function takes value $+\infty$ for every argument $x$. We will use Chung-Fuchs theorem to show under which condition subordinate random walk is transient. Denote with $\varphi_{X_1}$ the characteristic function of one step of a subordinate random walk. We want to prove that there exists $\delta > 0$ such that
\begin{equation*}
	\int_{{\langle -\delta, \delta \rangle}^d} {\re \OBL{\frac{1}{1 - \varphi_{X_1} (\theta)}} \dif \theta} < +\infty.
\end{equation*}
The law of variable $X_1$ is given with \eqref{eq:definition_of_xi}. We denote one step of the simple symmetric random walk $(Z_n)_{n \ge 0}$ with $Y_1$ and the characteristic function of that random variable with $\varphi$. Assuming $\aps{\theta} < 1$ we have
\begin{align*}
	\varphi_{X_1} (\theta)
	& = \bbE \UGL{e^{i \theta \cdot X_1}} = \sum_{x \in \bbZ^d} {e^{i \theta \cdot x} \sum_{m = 1}^{\infty} {\int_{\langle 0, +\infty \rangle} {\frac{t^m}{m!} e^{-t} \mu(\dif t)} \bbP(Z_m = x)}} \\
	& = \sum_{m = 1}^{\infty} {\int_{\langle 0, +\infty \rangle} {\frac{t^m}{m!} e^{-t} \mu(\dif t)} \sum_{x \in \bbZ^d} {e^{i \theta \cdot x} \bbP(Z_m = x)}} = \sum_{m = 1}^{\infty} {\int_{\langle 0, +\infty \rangle} {\frac{t^m}{m!} e^{-t} \mu(\dif t)} (\varphi(\theta))^m} \\
	& = \int_{\langle 0, +\infty \rangle} {\OBL{e^{t \varphi(\theta)} - 1} e^{-t} \mu(\dif t)} = \phi(1) - \phi(1 - \varphi(\theta)) = 1 - \phi(1 - \varphi(\theta)).
\end{align*}
From \cite[Section 1.2, page 13]{La96} we have
\begin{equation*}
	\varphi(\theta) = \frac{1}{d} \sum_{m = 1}^{d} {\cos(\theta_m)}, \quad \theta = (\theta_1, \theta_2, \ldots, \theta_m).
\end{equation*}
That is function with real values so
\begin{equation*}
	\int_{{\langle -\delta, \delta \rangle}^d} {\re \OBL{\frac{1}{1 - \varphi_{X_1}(\theta)}} \dif \theta} = \int_{{\langle -\delta, \delta \rangle}^d} {\frac{1}{\phi(1 - \varphi(\theta))} \dif \theta}.
\end{equation*}
From Taylor's theorem it follows that there exists $a \le 1$ such that
\begin{equation}\label{eq:taylor_ocjena_ch_f_od_Y}
	\aps{\varphi(\theta)} = \varphi(\theta) \le 1 - \frac{1}{4d} \aps{\theta}^2, \quad \theta \in B(0, a).
\end{equation}
Now we take $\delta$ such that ${\langle -\delta, \delta \rangle}^d \subset B(0, a)$. From \eqref{eq:taylor_ocjena_ch_f_od_Y}, using the fact that $\phi$ is increasing, we get
\begin{equation*}
	\frac{1}{\phi \OBL{1 - \varphi(\theta)}} \le \frac{1}{\phi \OBL{\aps{\theta}^2 / 4d}}, \quad \theta \in B(0, a).
\end{equation*}
Hence,
\begin{align*}
	\int_{{\langle -\delta, \delta \rangle}^d} {\frac{1}{\phi(1 - \varphi(\theta))} \dif \theta}
	& \le \int_{{\langle -\delta, \delta \rangle}^d} {\frac{1}{\phi \OBL{\aps{\theta}^2 / 4d}} \dif \theta} \le \int_{B(0, a)} {\frac{\phi(\aps{\theta}^2)}{\phi \OBL{\aps{\theta}^2 / 4d}} \frac{1}{\phi(\aps{\theta}^2)} \dif \theta} \\
	& \le a_2 (4d)^{\gamma_2} \int_{B(0, a)} {\frac{1}{\phi(\aps{\theta}^2)} \dif \theta} = c_1 (4d)^{\gamma_2} \int_{0}^{a} {\frac{r^{d - 1}}{\phi(r^2)} \dif r} \\
	& = \frac{c_1 (4d)^{\gamma_2}}{\phi(a)} \int_{0}^{a} {r^{d - 1} \frac{\phi(a)}{\phi(r^2)} \dif r} \le \frac{c_1 a_2 (4ad)^{\gamma_2}}{\phi(a)} \int_{0}^{a} {r^{d - 2 \gamma_2 - 1} \dif r}
\end{align*}
and the last integral converges for $d - 2 \gamma_2 - 1 > -1$. So, subordinate random walk is transient for $\gamma_2 < d/2$. Notice that in the case when $d \ge 3$ we have $\gamma_2 < d/2$ even if $\gamma_2 = 1$. That is the reason why we do not need \eqref{eq:upper_svojstvo_skaliranja} in proving results in dimensions higher than or equal to $3$. If some particular result depends on the dimension, we will write its proof using \eqref{eq:upper_svojstvo_skaliranja} just to show that the result is true even in dimensions $1$ and $2$ when $\gamma_2 < d/2$. In the case when $d \ge 3$, we can replace $a_2$ and $\gamma_2$ from \eqref{eq:upper_svojstvo_skaliranja} with $1$ and only use Lemma \ref{lm:free_upper_scaling}.

\subsection{Properties of functions $g$ and $j$}
\noindent
Let $\FJADEF{g}{\langle 0, +\infty \rangle}{\langle 0, +\infty \rangle}$ be defined by
\begin{equation}\label{eq:def_of_g}
	g(r) = \frac{1}{r^d \phi(r^{-2})}
\end{equation}
and let $\FJADEF{j}{\langle 0, +\infty \rangle}{\langle 0, +\infty \rangle}$ be defined by
\begin{equation}\label{eq:def_of_j}
	j(r) = r^{-d} \phi(r^{-2}).
\end{equation}
We use these functions in numerous places in our paper. In this section we present some of their properties that we need later.

\begin{LM}\label{lm:g_is_almost_decreasing}
	Assume \eqref{eq:upper_svojstvo_skaliranja} (if $d \le 2$) and let $1 \le r \le q$. Then $g(r) \ge a_2^{-1} g(q)$.
\end{LM}
\begin{PF}
	Using \eqref{eq:upper_svojstvo_skaliranja} and the fact that $d > 2\gamma_2$ we have
	\begin{equation*}
		g(r) = \frac{1}{\frac{r^d}{q^d} q^d \phi(q^{-2}) \frac{\phi(r^{-2})}{\phi(q^{-2})}} \ge \frac{1}{a_2} \OBL{\frac{q}{r}}^{d - 2\gamma_2} g(q) \ge \frac{1}{a_2} g(q).
	\end{equation*}
\end{PF}
We prove similar assertion for the function $j$.
\begin{LM}\label{lm:j_is_almost_decreasing}
	Assume \eqref{eq:svojstvo_skaliranja} and let $1 \le r \le q$. Then $j(r) \ge a_1 j(q)$.
\end{LM}
\begin{PF}
	Using \eqref{eq:svojstvo_skaliranja} we have
	\begin{equation*}
		j(r) = \frac{r^{-d}}{q^{-d}} q^{-d} \phi(q^{-2}) \frac{\phi(r^{-2})}{\phi(q^{-2})} \ge a_1 \OBL{\frac{q}{r}}^{d + 2\gamma_1} j(q) \ge a_1 j(q).
	\end{equation*}
\end{PF}
Using \eqref{eq:svojstvo_skaliranja}, \eqref{eq:upper_svojstvo_skaliranja} and Lemma \ref{lm:free_upper_scaling} we can easily prove a lot of different results about functions $g$ and $j$. We will state only those results that we need in the remaining part of our paper. For the first lemma we do not need any additional assumptions on the function $\phi$. For the second one we need \eqref{eq:svojstvo_skaliranja} and for the third one we need \eqref{eq:upper_svojstvo_skaliranja}.
\begin{LM}
	Let $r \ge 1$. If $0 < a \le 1$ then
	\begin{equation}\label{eq:j(ar)_le_noSC}
		j(ar) \le a^{-d - 2} j(r),
	\end{equation}
	\begin{equation}\label{eq:g(ar)_ge_noSC}
		g(ar) \ge a^{-d + 2} g(r).
	\end{equation}
	If $a \ge 1$ then
	\begin{equation}\label{eq:j(ar)_ge_noSC}
		j(ar) \ge a^{-d - 2} j(r).
	\end{equation}
\end{LM}

\begin{LM}
	Assume \eqref{eq:svojstvo_skaliranja} and let $0 < a \le 1$ and $r \ge 1$ such that $ar \ge 1$. Then
	\begin{equation}\label{eq:g(ar)_le_LSC}
		g(ar) \le \frac{g(r)}{a_1 a^{d - 2\gamma_1}}.
	\end{equation}
\end{LM}

\begin{LM}
	Assume \eqref{eq:upper_svojstvo_skaliranja} and let $r \ge 1$. If $0 < a \le 1$ such that $ar \ge 1$ then
	\begin{equation}\label{eq:g(ar)_ge_USC}
		g(ar) \ge \frac{g(r)}{a_2 a^{d - 2\gamma_2}}.
	\end{equation}		
	If $a \ge 1$ then
	\begin{equation}\label{eq:g(ar)_le_USC}
		g(ar) \le \frac{a_2}{a^{d - 2\gamma_2}} g(r).
	\end{equation}
\end{LM}

\section{Transition probability estimates}\label{sec:Transition_probability_estimates}
\noindent
In this section, we will investigate the behavior of the expression $\bbP(X_1 = z)$. We will prove that $\bbP(X_1 = z) \asymp j(\aps{z})$, $z \neq 0$. First we have to examine the behavior of the expression $c_m^{\phi}$.
\begin{LM}
	Assume \eqref{eq:svojstvo_skaliranja} and let $c_m^{\phi}$ be as in \eqref{eq:def_of_c_m^phi}. Then
		\begin{equation}\label{eq:asymptotics_of_c_m^phi}
		c_m^{\phi} \asymp \frac{\phi(m^{-1})}{m}, \quad m\in \bbN.
	\end{equation}
\end{LM}
\begin{PF}
	Since $\phi$ is a complete Bernstein function, there exists completely monotone density $\mu(t)$ such that
	\begin{equation*}
		c_m^{\phi} = \int_0^{\infty} {\frac{t^m}{m!} e^{-t} \mu(t) dt}, \quad m\in \bbN.
	\end{equation*}
	From \cite[Proposition 2.5]{KS14} we have
	\begin{equation}\label{eq:upper_bound_for_density_mu}
		\mu(t) \le (1 - 2e^{-1})^{-1} t^{-1} \phi(t^{-1}) = c_1 t^{-1} \phi(t^{-1}), \quad t > 0
	\end{equation}
	and
	\begin{equation}\label{eq:lower_bound_for_density_mu}
		\mu(t) \ge c_2 t^{-1} \phi(t^{-1}), \quad t \ge 1.
	\end{equation}
	Inequality \eqref{eq:lower_bound_for_density_mu} holds only if \eqref{eq:svojstvo_skaliranja} is satisfied and for inequality \eqref{eq:upper_bound_for_density_mu} we do not need any scaling conditions. Using monotonicity of $\mu$, \eqref{eq:incomplete_complete_gamma_ratio} and \eqref{eq:lower_bound_for_density_mu} we have
\begin{equation*}
	c_m^{\phi} \ge \frac{\mu(m)}{m!} \int_0^m {t^m e^{-t} dt} = \mu(m) \OBL{1 - \frac{\Gamma(m + 1, m)}{\Gamma(m + 1)}} \ge \frac{1}{4} \mu(m) \ge \frac{c_2}{4} \frac{\phi(m^{-1})}{m}
\end{equation*}
for $m$ large enough. On the other side,  using inequality \eqref{eq:upper_bound_for_density_mu}, monotonicity of $\mu$ and Lemma \ref{lm:free_upper_scaling}, we get
\begin{align*}
	c_m^{\phi}
	& \le \frac{1}{m!} \int_0^m {t^m e^{-t} c_1 \frac{\phi(t^{-1})}{t} dt} + \frac{\mu(m)}{m!} \int_m^{\infty} {t^m e^{-t} dt} \\
	& \le \frac{c_1}{m!} \phi(m^{-1}) \int_0^m {t^{m - 1} e^{-t} \frac{\phi(t^{-1})}{\phi(m^{-1})} dt} + \frac{\mu(m)}{m!} \int_0^{\infty} {t^m e^{-t} dt} \\
	& \le c_1 \phi(m^{-1}) \frac{1}{\Gamma(m)} \int_0^{\infty} {t^{m - 2} e^{-t} dt} + \mu(m) = c_1 \phi(m^{-1}) \frac{\Gamma(m - 1)}{\Gamma(m)} + \mu(m) \\
	& \le c_1 \frac{\phi(m^{-1})}{m} + c_1 \frac{\phi(m^{-1})}{m} = 2c_1 \frac{\phi(m^{-1})}{m}.
\end{align*}
Hence, we have
\begin{equation*}
	\frac{c_2}{4} \frac{\phi(m^{-1})}{m} \le c_m^{\phi} \le 2c_1 \frac{\phi(m^{-1})}{m}
\end{equation*}
for $m$ large enough, but we can change constants and get \eqref{eq:asymptotics_of_c_m^phi}.
\end{PF}
We are now ready to examine the expression $\bbP(X_1 = z)$.
\begin{PROP}\label{prop:asymptotics_for_jumps}
	Assume \eqref{eq:svojstvo_skaliranja}. Then
	\begin{equation*}
		\bbP(X_1 = z) \asymp \aps{z}^{-d} \phi(\aps{z}^{-2}), \quad z\neq 0.
	\end{equation*}
\end{PROP}
\begin{PF}
Using \eqref{eq:definition_of_xi} and the fact that $\bbP(Z_m = z) = 0$ for $\aps{z} > m$, we have
\begin{equation*}
	\bbP(X_1 = z) = \sum_{m \ge \aps{z}} {c_m^{\phi} \bbP(Z_m = z)}.
\end{equation*}
To get the upper bound for the expression $\bbP(X_1 = z)$ we will use \cite[Theorem 2.1]{HS93} which states that there are constants $C' > 0$ and $C > 0$ such that
\begin{equation}\label{eq:gaussian_upper_bound_from_HSC}
	\bbP(Z_m = z) \le C' m^{-\frac{d}{2}} e^{-\frac{\aps{z}^2}{C m}}, \quad \forall\, z \in \bbZ^d,\, \forall\, m \in \bbN.
\end{equation}
Together with \eqref{eq:asymptotics_of_c_m^phi} this result yields
\begin{align*}
	\bbP(X_1 = z)
	& \le \sum_{m \ge \aps{z}} {c_1 \frac{\phi(m^{-1})}{m} C' m^{-\frac{d}{2}} e^{-\frac{\aps{z}^2}{C m}}} \le c_2 \int_{\aps{z}}^{\infty} {\phi(t^{-1}) t^{-\frac{d}{2} - 1} e^{-\frac{\aps{z}^2}{C t}} dt} \\
	& = c_2 \int_0^{\frac{\aps{z}}{C}} {\phi(C s\aps{z}^{-2}) \OBL{\frac{\aps{z}^2}{C s}}^{-\frac{d}{2} - 1} e^{-s} \ \frac{\aps{z}^2}{C s^2} ds} \\
	& = c_3 \aps{z}^{-d} \OBL{\int_0^{\frac{1}{C}} {\phi(C s\aps{z}^{-2}) s^{\frac{d}{2} - 1} e^{-s} ds} + \int_{\frac{1}{C}}^{\frac{\aps{z}}{C}} {\phi(C s\aps{z}^{-2}) s^{\frac{d}{2} - 1} e^{-s} ds}} \\
	& =: c_3 \aps{z}^{-d} (I_1(z) + I_2(z)).
\end{align*}
Let us first examine $I_1(z)$. Using \eqref{eq:svojstvo_skaliranja}, we get
\begin{equation*}
	I_1(z) = \phi(\aps{z}^{-2}) \int_0^{\frac{1}{C}} {\frac{\phi (C s\aps{z}^{-2})}{\phi(\aps{z}^{-2})} s^{\frac{d}{2} - 1} e^{-s} ds} \le c_4 \phi(\aps{z}^{-2}).
\end{equation*}
Using Lemma \ref{lm:free_upper_scaling} instead of \eqref{eq:svojstvo_skaliranja} and replacing the upper limit in the integral by $\infty$, we get $I_2(z) \le c_5 \phi(\aps{z}^{-2})$. Hence, $\bbP(X_1 = z) \le c_6 \aps{z}^{-d} \phi(\aps{z}^{-2})$.

In finding the matching lower bound for $\bbP(X_1 = z)$, periodicity of a simple random walk plays very important role. We write $n \leftrightarrow x$ if $n$ and $x$ have the same parity, i.e., if $n + x_1 + x_2 + \cdots + x_d$ is even. Directly from \cite[Proposition 1.2.5]{La96}, we get
\begin{equation}\label{eq:lower_bound_A}
	\bbP(Z_m = z) \ge c_7 m^{-\frac{d}{2}} e^{-\frac{d \aps{z}^2}{2m}}
\end{equation}
for $0 \leftrightarrow z \leftrightarrow m$ and $\aps{z} \le m^{\alpha}$, $\alpha < 2/3$. In the case when $1 \leftrightarrow z \leftrightarrow m$ we have
\begin{equation}\label{eq:pom_for_z_odd}
	\bbP(Z_m = z) = \frac{1}{2d} \sum_{i = 1}^d {[\bbP(Z_{m - 1} = z + e_i) + \bbP(Z_{m - 1} = z - e_i)]}.
\end{equation}
By combining \eqref{eq:lower_bound_A} and \eqref{eq:pom_for_z_odd}, we can easily get
\begin{equation}\label{eq:lower_bound_B}
	\bbP(Z_m = z) \ge c_{8} m^{-\frac{d}{2}} e^{-\frac{\aps{z}^2}{c m}}, \quad \aps{z} \le m^{\frac{1}{2}}, 1 \leftrightarrow z \leftrightarrow m.
\end{equation}
We will find lower bound for $\bbP(X_1 = z)$ when $z \leftrightarrow 0$  by using \eqref{eq:lower_bound_A}, the proof when $z \leftrightarrow 1$ being analogous using \eqref{eq:lower_bound_B}. If $z \leftrightarrow 0$ then $\bbP(Z_m = z) =  0$ for $m = 2l - 1$, $l \in \bbN$. Hence,
\begin{align*}
	\bbP(X_1 = z)
	& \ge \sum_{m \ge \aps{z}^2, m = 2l} {c_{9} \frac{\phi(m^{-1})}{m}} m^{-\frac{d}{2}} e^{-\frac{-d\aps{z}^2}{2m}} = c_{9} \sum_{l \ge \aps{z}^2 / 2} {\frac{\phi((2l)^{-1})}{2l} (2l)^{-\frac{d}{2}} e^{-\frac{d\aps{z}^2}{4l}}} \\
	& \ge c_{10} \int_{\aps{z}^2 / 2}^{\infty} {\frac{\phi((2t)^{-1})}{2t} (2t)^{-\frac{d}{2}} e^{-\frac{d\aps{z}^2}{4t}} \dif t} = \frac{c_{10}}{2} \int_{\aps{z}^2}^{\infty} {\phi(t^{-1}) t^{-\frac{d}{2} - 1} e^{-\frac{d\aps{z}^2}{2t}} \dif t} \\
	& = \frac{c_{10}}{2} \int_0^{\frac{d}{2}} {\phi \OBL{\frac{2s}{d \aps{z}^2}} \OBL{\frac{d \aps{z}^2}{2s}}^{-\frac{d}{2} - 1} e^{-s}\ \frac{d \aps{z}^2}{2s^2} \dif s} \\
	& = c_{11} \aps{z}^{-d} \phi(\aps{z}^{-2}) \int_0^{\frac{d}{2}} {\frac{\phi \OBL{\frac{2s}{d} \aps{z}^{-2}}}{\phi (\aps{z}^{-2})} s^{\frac{d}{2} - 1} e^{-s} \dif s} \\
	& \ge c_{11} \aps{z}^{-d} \phi(\aps{z}^{-2}) \int_0^{\frac{d}{2}} {\frac{2}{d} s^{\frac{d}{2}} e^{-s} ds} = c_{12} \aps{z}^{-d} \phi(\aps{z}^{-2}),
\end{align*}
where in the last line we used Lemma \ref{lm:free_upper_scaling}.
\end{PF}
\begin{REM}\label{rem:infinite_variance_of_the_step}
It follows immediately form Proposition \ref{prop:asymptotics_for_jumps} that the second moment of the step $X_1$ is infinite.
\end{REM}

\section{Green function estimates}\label{sec:Green_function_estimates}
\noindent
The Green function of $X$ is defined by $G(x, y) = G(y - x)$, where
\begin{equation*}
	G(y) = \bbE \left[\sum\limits_{n = 0}^{\infty} {\bbjedan_{\{X_n = y\}}}\right].
\end{equation*}
Note that for $n \ge 1$
\begin{align*}
	\bbP(X_n = y)
		& = \bbP(Z_{T_n} = y) = \sum_{m = n}^{\infty} {\bbP(Z_m = y) \bbP(T_n = m)} \\
		& = \sum_{m = n}^{\infty} {\sum_{m_1 + \cdots + m_n = m} {c_{m_1}^{\phi} \cdots c_{m_n}^{\phi}} \bbP(Z_m = y)}
	\end{align*}
Hence, for $y \neq 0$ we have
\begin{equation}\label{eq:expression_for_G}
	G(y) = \sum\limits_{m = 1}^{\infty} {\sum\limits_{n = 1}^{m} \sum\limits_{m_1 + \cdots + m_n = m} {c_{m_1}^{\phi} \cdots c_{m_n}^{\phi} \bbP(Z_m = y)}} = \sum\limits_{m = 1}^{\infty} {c(m) \bbP(Z_m = y)},
\end{equation}
where
\begin{equation}\label{eq:definition_of_c(m)}
	c(m) = \sum\limits_{n = 1}^{m} \sum\limits_{m_1 + \cdots + m_n = m} {c_{m_1}^{\phi} \cdots c_{m_n}^{\phi}} = \sum_{n = 0}^{\infty} {\bbP(T_n = m)},
\end{equation}
and $T_n$ is as before. Now we will investigate the behavior of the sequence $c(m)$. Since $\phi$ is a complete Bernstein function (hence special), we have
	\begin{equation}\label{eq:introduction_of_u}
		\frac{1}{\phi(\lambda)} = \int_{\langle 0, \infty \rangle} {e^{-\lambda t} u(t) \dif t}
	\end{equation}
	for some non-increasing function $\FJADEF{u}{\langle 0, \infty \rangle}{\langle 0, \infty \rangle}$ satisfying $\int_0^1{u(t)\dif t} < \infty$, see \cite[Theorem 11.3.]{SS12}.

\begin{LM}\label{lm:expression_for_c(m)}
	Let $c(m)$ be as in \eqref{eq:definition_of_c(m)}. Then
\begin{equation}\label{eq:expression_for_c(k)}
		c(m) = \frac{1}{m!} \int_{\langle 0, \infty \rangle} {t^m e^{-t} u(t) \dif t}, \quad m \in \bbN_0.
	\end{equation}
\end{LM}
\begin{PF}
	We follow the proof of \cite[Theorem 2.3]{BC10}. Define $M(x) = \sum_{m \le x} {c(m)}$, $x \in \bbR$. The Laplace transformation $\calL(M)$ of the measure generated by $M$ is equal to	\begin{align}\label{al:LM_lambda_prvi_nacin}
		\calL(M)(\lambda)
		& = \int_{[0, \infty\rangle} {e^{-\lambda x} \dif M(x)} = \sum_{m = 0}^{\infty} {c(m) e^{-\lambda m}} = \sum_{m = 0}^{\infty} {e^{-\lambda m}} \sum_{n = 0}^{\infty} {\bbP(T_n = m)} \nonumber \\
		& = \sum_{n = 0}^{\infty} \sum_{m = 0}^{\infty} {e^{-\lambda m} \bbP(T_n = m)} = \sum_{n = 0}^{\infty} {\bbE[e^{-\lambda T_n}]} = \sum_{n = 0}^{\infty} {\OBL{\bbE[e^{-\lambda R_1}]}^n} = \frac{1}{1 - \bbE[e^{-\lambda R_1}]}.
	\end{align}
	Now we calculate $\bbE[e^{-\lambda R_1}]$:
	\begin{align*}
		\bbE[e^{-\lambda R_1}]
		& = \sum_{m = 1}^{\infty} {e^{-\lambda m} \int_{\langle 0, \infty \rangle} {\frac{t^m}{m!} e^{-t}} \mu(\dif t)} = \int_{\langle 0, \infty \rangle} {\OBL{e^{te^{-\lambda}} - 1} e^{-t} \mu(\dif t)} = 1 - \phi(1 - e^{-\lambda}),
	\end{align*}
	where we used $\phi(1) = 1$ in the last equality. Hence, $\calL(M)(\lambda) = 1 / \phi(1 - e^{-\lambda})$. Now we define $\Phi(\lambda) := 1 / \phi(\lambda)$ and we want to show that
	\begin{equation*}
		\Phi(1 - e^{-\lambda}) = \sum_{m = 0}^{\infty} {\frac{(-1)^m \Phi^{(m)}(1)}{m!} e^{-\lambda m}}.
	\end{equation*}
	It is easy to see that
	$\Phi^{(m)}(\lambda) = (-1)^m \int_{\langle 0, \infty \rangle} {t^m e^{-\lambda t}u(t) \dif t}$. Hence,	\begin{equation}\label{al:LM_lambda_drugi_nacin}
		\sum_{m = 0}^{\infty} {\frac{(-1)^m \Phi^{(m)}(1)}{m!} e^{-\lambda m}} = \sum_{m = 0}^{\infty} {\frac{(-1)^m}{m!} (-1)^m \int_{\langle 0, \infty \rangle} {t^m e^{-t} u(t) \dif t} e^{-\lambda m}} = \Phi(1 - e^{-\lambda}).
	\end{equation}
	Since $\calL(M)(\lambda) = 1/\phi(1 - e^{-\lambda}) = \Phi(1 - e^{-\lambda})$ from calculations \eqref{al:LM_lambda_prvi_nacin} and \eqref{al:LM_lambda_drugi_nacin} we have
	\begin{equation*}
		\sum_{m = 0}^{\infty} {c(m) e^{-\lambda m}} = \sum_{m = 0}^{\infty} {\frac{1}{m!} \int_{\langle 0, \infty \rangle} {t^m e^{-t} u(t) \dif t} \, e^{-\lambda m}}.
	\end{equation*}
	The statement of this lemma follows by the uniqueness of the Laplace transformation.
\end{PF}

\begin{LM}\label{lm:asimptotika_od_c(m)}
	Assume \eqref{eq:svojstvo_skaliranja}. Then
	\begin{equation*}
		c(m) \asymp \frac{1}{m \phi(m^{-1})}, \quad m \in \bbN.
	\end{equation*}
\end{LM}
\begin{PF}
	Since $\phi$ is a complete Bernstein function, using \eqref{eq:svojstvo_skaliranja} we can obtain
		\begin{equation}\label{eq:u(t)_comparable}
		u(t) \asymp \frac{1}{t \phi(t^{-1})}, \quad t \ge 1,
	\end{equation}
	where the upper bound is valid even without \eqref{eq:svojstvo_skaliranja} (see \cite[Corollary 2.4.]{KS14}) and $u$ is as in Lemma \ref{lm:expression_for_c(m)}. Using monotonicity of $u$, \eqref{eq:incomplete_complete_gamma_ratio} and \eqref{eq:u(t)_comparable}, we get that
	\begin{equation*}
	c(m) \ge u(m) \frac{1}{m!} \int_0^m {t^m e^{-t} \dif t} = u(m) \OBL{1 - \frac{\Gamma(m + 1, m)}{\Gamma(m + 1)}} \ge \frac{1}{4} u(m) \ge \frac{c_1}{m \phi(m^{-1})},
	\end{equation*}
	for $m$ large enough. Now we will find the upper bound for $c(m)$. Here we use that $t \mapsto t^m e^{-t}$ is unimodal with maximum at $m$. By splitting the integral and using \eqref{eq:stirlings_formula}, we have
	\begin{align*}
		c(m)
		& = \frac{1}{m!} \int_{0}^{m/2} {t^m e^{-t} u(t) \dif t} + \frac{1}{m!} \int_{m/2}^{m} {t^m e^{-t} u(t) \dif t} + \frac{1}{m!} \int_{m}^{\infty} {t^m e^{-t} u(t) \dif t} \\
		& \le \frac{c_2 2^{-m} m^m e^{-m/2}}{\sqrt{2 \pi m} m^m e^{-m}} \int_0^{m/2} {u(t) \dif t} + u(m/2) \frac{1}{m!} \int_{m/2}^{\infty} {t^m e^{-t} \dif t} + u(m) \frac{1}{m!} \int_m^{\infty} {t^m e^{-t} \dif t} \\ 
		& \le c_2 \frac{(2^{-1} e^{1/2})^m}{\sqrt{2 \pi m}} \int_0^1{u(t) \dif t} + c_2 \frac{(2^{-1} e^{1/2})^m}{\sqrt{2 \pi m}} \int_1^{m/2} {u(t) \dif t} + u(m/2) + u(m),
	\end{align*}
	for $m$ large enough. Since $2^{-1} e^{1/2} < 1$, $\int_0^1{u(t) \dif t} < \infty$ and $\phi$ is increasing, we have
	\begin{equation*}
		c_2 \frac{(2^{-1} e^{1/2})^m}{\sqrt{2 \pi m}} \int_0^1{u(t) \dif t} \le \frac{1}{m} \le \frac{1}{m \phi(m^{-1})} \le \frac{1}{4c_1} u(m)
	\end{equation*}
	for $m$ large enough, where we used \eqref{eq:u(t)_comparable} in the last inequality. We will estimate the integral $\int_1^{m/2} {u(t) \dif t}$ by $mu(m/2)$. Using \eqref{eq:u(t)_comparable} and \eqref{eq:svojstvo_skaliranja} we get
	\begin{align*}
		\int_1^{m/2} {u(t)} dt
		& \le c_3 \int_1^{m/2} {\frac{1}{t\phi(t^{-1})}} dt = \frac{c_3}{\phi(2m^{-1})} \int_1^{m/2} {\frac{\phi(2m^{-1})}{\phi(t^{-1})} t^{-1}} dt \\
		& \le \frac{c_3}{\phi(2m^{-1})} \frac{1}{a_1(m/2)^{\gamma_1}} \int_1^{m/2} {t^{\gamma_1 - 1}} dt \\
		& \le \frac{c_3}{a_1\phi(2m^{-1})(m/2)^{\gamma_1}} \frac{(m/2)^{\gamma_1}}{\gamma_1} = \frac{c_3}{a_1\gamma_1} \frac{1}{\phi(2m^{-1})}.
	\end{align*}
	Since $u(m/2) \ge c_4 / ((m/2)\phi(2m^{-1}))$ for $m \ge 2$, we have $1 / \phi(2m^{-1}) \le (1 / 2c_4) mu(m/2)$. Hence,
	\begin{equation*}
		\int_1^{m/2} {u(t)} dt \le \frac{c_3}{a_1\gamma_1} \frac{1}{2c_4} mu(m/2) = c_5mu(m/2).
	\end{equation*}
	Using this estimate and the fact that $2^{-1} e^{1/2}$ is less then $1$, we have
	\begin{equation*}
		c_2 \frac{(2^{-1} e^{1/2})^m}{\sqrt{2\pi m}} \int_1^{m/2} {u(t)} dt \le c_6 \frac{(2^{-1} e^{1/2})^m}{\sqrt{2\pi m}} m u(m/2) = c_6(2^{-1} e^{1/2})^m m^{1/2} u(m/2) \le u(m/2)
	\end{equation*}
	for $m$ large enough. Now we have to show that $u(m/2)$ can be estimated by $u(m)$. Again, we will use \eqref{eq:u(t)_comparable} and \eqref{eq:svojstvo_skaliranja}:
	\begin{equation*}
		u(t/2) \le \frac{c_3}{(t/2)\phi(2t^{-1})}  = \frac{2c_3}{t} \frac{\phi(t^{-1})}{\phi(2t^{-1})} \frac{1}{\phi(t^{-1})} \le \frac{2c_3}{a_1 2^{\gamma_1}} \frac{1}{t \phi(t^{-1})} \le \frac{2c_3}{c_4 a_1 2^{\gamma_1}} u(t) = c_7 u(t),
	\end{equation*}
	where we assumed that $t \ge 2$ because we need $2t^{-1} \le 1$ so that we can use \eqref{eq:svojstvo_skaliranja}. Now, for $m$ large enough, we have
	\begin{equation*}
		c(m) \le \frac{1}{4c_1} u(m) + u(m/2) + u(m/2) + u(m) \le \frac{1}{4c_1} u(m) + 2c_7 u(m) + u(m) \le \frac{c_8}{m \phi(m^{-1})}.
	\end{equation*}
	Hence,
	\begin{equation*}
		\frac{c_1}{m \phi(m^{-1})} \le c(m) \le \frac{c_8}{m \phi(m^{-1})}
	\end{equation*}
	for $m$ large enough. We can now change constants in such a way that the statement of this lemma is true for every $m \in \bbN$.
\end{PF}
\begin{TM}\label{tm:asymptotics_of_G}
	Assume \eqref{eq:svojstvo_skaliranja} and, if $d \le 2$, assume additionally \eqref{eq:upper_svojstvo_skaliranja}. Then
	\begin{equation}\label{eq:asimptotika_od_G}
		G(x) \asymp \frac{1}{\aps{x}^d \phi(\aps{x}^{-2})}, \quad \aps{x} \ge 1.
	\end{equation}
\end{TM}
\begin{PF}
	We assume $\aps{x} \ge 1$ throughout the whole proof. In \eqref{eq:expression_for_G} we showed that $G(x) = \sum_{m = 1}^{\infty} {c(m) p(m, x)}$, where $p(m, x) = \bbP(Z_m = x)$. Let $q(m, x) = 2 \OBL{d / (2 \pi m)}^{\frac{d}{2}} e^{-\frac{d \aps{x}^2}{2m}}$ and $E(m, x) = p(m, x) - q(m, x)$. By \cite[Theorem 1.2.1]{La96}
	\begin{equation}\label{eq:upper_bound_for_E}
		\aps{E(m, x)} \le c_1 m^{-\frac{d}{2}} / \aps{x}^2.
	\end{equation}
	Since $p(m, x) = 0$ for $m <\aps{x}$, we have
	\begin{equation*}
		G(x) = \sum_{m > \aps{x}^2} {c(m) p(m, x)} + \sum_{\aps{x} \le m \le \aps{x}^2} {c(m) p(m, x)} =: J_1(x) + J_2(x).
	\end{equation*}
	First we estimate
	\begin{equation*}
		J_1(x) = \sum_{m > \aps{x}^2} {c(m) q(m, x)} + \sum_{m > \aps{x}^2} {c(m) E(m, x)} =: J_{11}(x) + J_{12}(x).
	\end{equation*}
	By Lemma \ref{lm:asimptotika_od_c(m)}, \eqref{eq:upper_bound_for_E} and \eqref{eq:upper_svojstvo_skaliranja}
	\begin{align*}
		\aps{J_{12}(x)}
		& \le c_2 \sum_{m > \aps{x}^2} {\frac{1}{m \phi(m^{-1})} \frac{m^{-\frac{d}{2}}}{\aps{x}^2}} = \frac{c_2}{\aps{x}^2 \phi(\aps{x}^{-2})} \sum_{m > \aps{x}^2} {\frac{\phi(\aps{x}^{-2})}{\phi(m^{-1})} m^{-\frac{d}{2} - 1}} \\
		& \le \frac{c_3 \aps{x}^{-2 \gamma_2}}{\aps{x}^2 \phi(\aps{x}^{-2})} \int_{\aps{x}^2}^{\infty} {t^{\gamma_2 - \frac{d}{2} - 1} \dif t} = \frac{c_4}{\aps{x}^2} \frac{1}{\aps{x}^d \phi(\aps{x}^{-2})}.
	\end{align*}
	Now we have
	\begin{equation*}
		\lim_{\aps{x} \to \infty} {\aps{x}^d \phi(\aps{x}^{-2}) \aps{J_{12}(x)}} = 0.
	\end{equation*}
	By Lemma \ref{lm:asimptotika_od_c(m)}, \eqref{eq:svojstvo_skaliranja} and \eqref{eq:upper_svojstvo_skaliranja}
	\begin{align*}
		J_{11}(x)
		& \asymp \int_{\aps{x}^2}^{\infty} {\frac{1}{t \phi(t^{-1})} t^{-\frac{d}{2}} e^{-\frac{d \aps{x}^2}{2 t}} \dif t} = \frac{1}{\phi(\aps{x}^{-2})} \int_{\aps{x}^2}^{\infty} {\frac{\phi(\aps{x}^{-2})}{\phi(t^{-1})} t^{-\frac{d}{2} - 1} e^{-\frac{d \aps{x}^2}{2 t}} \dif t} \\
		& \asymp \frac{\aps{x}^{-2 \gamma_i}}{\phi(\aps{x}^{-2})} \int_{\aps{x}^2}^{\infty} {t^{\gamma_i - \frac{d}{2} - 1} e^{-\frac{d \aps{x}^2}{2 t}} \dif t} = \frac{1}{\aps{x}^d \phi(\aps{x}^{-2})} \int_0^{\frac{d}{2}} {s^{\frac{d}{2} -\gamma_i - 1} e^{-s} \dif s} \asymp \frac{1}{\aps{x}^d \phi(\aps{x}^{-2})},
	\end{align*}
	where the last integral converges because of the condition $\gamma_2 < d/2$. We estimate $J_2(x)$ using \eqref{eq:gaussian_upper_bound_from_HSC} and \eqref{eq:svojstvo_skaliranja}:
	\begin{align*}
		J_2(x)
		& \le c_5  \int_{\aps{x}}^{\aps{x}^2} {\frac{t^{-\frac{d}{2} - 1}}{\phi(t^{-1})} e^{-\frac{\aps{x}^2}{C t}} \dif t} = \frac{c_5}{\phi(\aps{x}^{-2})} \int_{{\aps{x}}}^{\aps{x}^2} {\frac{\phi(\aps{x}^{-2})}{\phi(t^{-1})} t^{-\frac{d}{2} - 1} e^{-\frac{\aps{x}^2}{C t}} \dif t} \\
		& \le \frac{c_5 \aps{x}^{-2 \gamma_1}}{a_1 \phi(\aps{x}^{-2})} \int_{{\aps{x}}}^{\aps{x}^2} {t^{\gamma_1 - \frac{d}{2} - 1} e^{-\frac{\aps{x}^2}{C t}} \dif t} = \frac{c_5 \aps{x}^{-2 \gamma_1}}{a_1 \phi(\aps{x}^{-2})} \int_{\frac{1}{C}}^{\frac{\aps{x}}{C}} {\OBL{\frac{\aps{x}^2}{Cs}}^{\gamma_1 - \frac{d}{2} - 1} e^{-s} \frac{\aps{x}^2}{C s^2} \dif s} \\
		& \le \frac{c_6}{\aps{x}^d \phi(\aps{x}^{-2})} \int_0^{\infty} {s^{\frac{d}{2} - \gamma_1 - 1} e^{-s} \dif s} = \frac{c_7}{\aps{x}^d \phi(\aps{x}^{-2})}.
	\end{align*}
	Using $J_{11}(x) \ge (2 c_8) / (\aps{x}^{d} \phi(\aps{x}^{-2}))$ and $J_{12}(x) \aps{x}^d \phi(\aps{x}^{-2}) \ge - c_8$ for $\aps{x}$ large enough and for some constant $c_8 > 0$, we get
	\begin{equation*}
		G(x) \aps{x}^d \phi(\aps{x}^{-2}) \ge J_{11}(x) \aps{x}^d \phi(\aps{x}^{-2}) + J_{12}(x) \aps{x}^d \phi(\aps{x}^{-2}) \ge 2c_8 - c_8 = c_8
	\end{equation*}
	On the other hand
	\begin{equation*}
		G(x) \aps{x}^d \phi(\aps{x}^{-2}) \le c_{9} +  J_{12}(x) \aps{x}^d \phi(\aps{x}^{-2}) + c_7 \le 2 c_{9} + c_7 = c_{10}.
	\end{equation*}
	Here we used $J_{11}(x) \le c_9 / (\aps{x}^d \phi(\aps{x}^{-2}))$, $J_2(x) \le c_7 / (\aps{x}^d \phi(\aps{x}^{-2}))$ and $J_{12}(x) \aps{x}^d \phi(\aps{x}^{-2}) \\ \le c_{9}$ for $\aps{x}$ large enough and for some constant $c_9 > 0$. So, we have $c_8 \le G(x) \aps{x}^d \phi(\aps{x}^{-2}) \le c_{10}$ for $\aps{x}$ large enough. Now we can change the constants $c_8$ and $c_{10}$ to get	\begin{equation*}
		G(x) \asymp \frac{1}{\aps{x}^d \phi(\aps{x}^{-2})}, \quad \textnormal{for all}\,\, \aps{x} \ge 1.
	\end{equation*}
\end{PF}

\section{Estimates of the Green function of a ball}\label{sec:Estimates_of_Green_function_of_a_ball}
\noindent
Let $B \subset \bbZ^d$ and define
\begin{equation*}
	G_B (x, y) = \bbE_x \UGL{\sum\limits_{n = 0}^{\tau_B - 1} {\bbjedan_{\{X_n = y\}}}}
\end{equation*}
where $\tau_B$ is as before. A well-known result about Green function of a set is formulated in the following lemma.
\begin{LM}\label{lm:Green_function_of_a_ball}
	Let $B$ be a finite subset of $\bbZ^d$. Then
	\begin{align*}
		G_B (x, y) & = G(x, y) - \bbE_x \UGL{G(X_{\tau_B}, y)}, \quad x, y \in B, \\
		G_B (x, x) & = \frac{1}{\bbP_x (\tau_B < \sigma_x)}, \quad x\in B,		
	\end{align*}
	where $\sigma_x = \inf \{n \ge 1 : X_n = x\}$.
\end{LM}
Our approach in obtaining estimates for the Green function of a ball uses the maximum principle for the operator $A$ that we define by
\begin{equation}\label{eq:def_od_A}
	(Af)(x) := ((P - I)f)(x) = (Pf)(x) - (If)(x) = \sum_{y \in \bbZ^d} {p(x, y) f(y)} - f(x).
\end{equation}
Since $\sum_{y \in \bbZ^d} {p(x, y)} = 1$ and $p(x, y) = \bbP(X_1 = y - x)$ we have
\begin{equation*}
	(Af)(x) = \sum_{y \in \bbZ^d} {\bbP(X_1 = y - x) (f(y) - f(x))}.
\end{equation*}
Before proving the maximum principle, we will show that for the function $\eta(x) := \bbE_x [\tau_{B_n}]$ we have $(A\eta)(x) = -1$, for all $x \in B_n$. Let $x \in B_n$. Then
\begin{equation*}
	\eta(x) = \sum_{y \in \bbZ^d} {\bbE_x[\tau_{B_n} \mid X_1 = y] \bbP_x(X_1 = y)} = \sum_{y \in \bbZ^d} {(1 + \bbE_y[\tau_{B_n}]) \bbP(X_1 = y - x)} = 1 + (P\eta)(x)
\end{equation*}
and this is obviously equivalent to $(A\eta)(x) = -1$, for all $ x \in B_n$. It follows from the Definition \ref{def:harmonic_function} that $f$ is harmonic in $B \subset \bbZ^d$ if and only if $(Af)(x) = 0$, for all $x \in B$.

\begin{PROP}\label{tm:maximum_principle} Assume that there exists $x \in \bbZ^d$ such that $(Af)(x) < 0$. Then
	\begin{equation}\label{eq:max_principle}
		f(x) > \inf_{y \in \bbZ^d} f(y).
	\end{equation}		
\end{PROP}
\begin{PF}
	If \eqref{eq:max_principle} is not true, then $f(x) \le f(y)$, for all $ y \in \bbZ^d$. In this case, we have
	\begin{equation*}
		(Pf)(x) = \sum_{y \in \bbZ^d} {\bbP(X_1 = y - x) f(y)} \ge f(x) \sum_{y \in \bbZ^d} {\bbP(X_1 = y - x)} = f(x).
	\end{equation*}
	This implies $(Af)(x) = (Pf)(x) - f(x) \ge 0$ which is in contradiction with the assumption that $(Af)(x) < 0$.
\end{PF}
We will now prove a series of lemmas and propositions in order to get the estimates for the Green function of a ball. In all those results we assume \eqref{eq:svojstvo_skaliranja} and, if $d \le 2$, we additionally assume \eqref{eq:upper_svojstvo_skaliranja}. Throughout the rest of this section, we follow \cite[Section 4]{MK12}.
\begin{LM}\label{lm:G_(B_n)>=G}
	There exist $a \in \langle 0, 1/3 \rangle$ and $C_1 > 0$ such that for every $n \in \bbN$
	\begin{equation}\label{eq:G_(B_n)>=G}
		G_{B_n}(x, y) \ge C_1 G(x, y), \quad \forall\, x, y \in B_{an}.
	\end{equation}
\end{LM}
\begin{PF}
	From Lemma \ref{lm:Green_function_of_a_ball} we have
	\begin{equation*}
		G_{B_n} (x, y) = G(x, y) - \bbE_x[G(X_{\tau_{B_n}}, y)].
	\end{equation*}
	We will first prove this lemma in the case when $x \neq y$. If we show that $\bbE_x[G(X_{\tau_{B_n}}, y)] \le c_1 G(x, y)$ for some $c_1 \in \langle 0, 1 \rangle$ we will have \eqref{eq:G_(B_n)>=G} with the constant $c_2 = 1 - c_1$. Let $a \in \langle 0, 1/3 \rangle$ and $x, y \in B_{an}$. In that case, we have $\aps{x - y} \le 2an$. Since $X_{\tau_{B_n}} \notin B_n$, $x \neq y$ and $(1 - a) / (2a) > 1$ if and only if $a < 1/3$, we have
	\begin{equation}\label{eq:aps(y-X_tau_Bn)}
		\aps{y - X_{\tau_{B_n}}} \ge (1 - a)n = \frac{1 - a}{2a} 2an \ge \frac{1 - a}{2a} \aps{x - y} \ge 1.
	\end{equation}
	Using Theorem \ref{tm:asymptotics_of_G}, \eqref{eq:aps(y-X_tau_Bn)}, Lemma \ref{lm:g_is_almost_decreasing} and \eqref{eq:g(ar)_le_USC}, we get
	\begin{align*}
		G(X_{\tau_{B_n}}, y)
		& \asymp g(\aps{y - X_{\tau_{B_n}}}) \le a_2 g\OBL{\frac{1 - a}{2a} \aps{x - y}} \\
		& \le a_2^2 \OBL{\frac{2a}{1 - a}}^{d - 2\gamma_2} g(\aps{x - y}) \asymp a_2^2 \OBL{\frac{2a}{1 - a}}^{d - 2\gamma_2} G(x, y).
	\end{align*}
	Since $2a / (1 - a) \longrightarrow 0$ when $a \rightarrow 0$ and $d > 2 \gamma_2$, if we take $a$ small enough and then fix it, we have $\bbE_x[G(X_{\tau_{B_n}}, y)] \le c_1 G(x, y)$ for $c_1 \in \langle 0, 1 \rangle$ and that is what we wanted to prove. Now we deal with the case when $x = y$. From Lemma \ref{lm:Green_function_of_a_ball} we have $G_{B_n}(x, x) = (\bbP(\tau_{B_n} < \sigma_x))^{-1}$ and from the definition of the function $G$ and the transience of random walk we get $G(x, x) = G(0) \in [1, \infty \rangle$. Now, we can conclude that
	\begin{equation*}
		G_{B_n}(x, x) \ge 1 = (G(0))^{-1} G(0) = (G(0))^{-1} G(x, x).
	\end{equation*}
	If we define $C_1 := \min\{c_2, (G(0))^{-1}\}$ we have \eqref{eq:G_(B_n)>=G}.
\end{PF}
Using Lemma \ref{lm:G_(B_n)>=G} we can prove the following result:
\begin{PROP}\label{prop:E[tau_Bn]>=}
	There exists constant $C_2 > 0$ such that for all $n \in \bbN$
	\begin{equation}\label{eq:E[tau_Bn]>=}
		\bbE_x[\tau_{B_n}] \ge \frac{C_2}{\phi(n^{-2})}, \quad \forall\, x \in B_{\frac{an}{2}},
	\end{equation}
	where $a \in \langle 0, 1/3 \rangle$ is as in Lemma \ref{lm:G_(B_n)>=G}.
\end{PROP}
\begin{PF}
	Let $x \in B_{\frac{an}{2}}$. In that case, we have $B(x, an/2) \subseteq B_{an}$. We set $b = a/2$ for easier notation. Notice that $\bbE_x[\tau_{B_n}] = \sum_{y \in B_n} G_{B_n}(x, y)$. Using this equality, Lemma \ref{lm:G_(B_n)>=G}, Theorem \ref{tm:asymptotics_of_G} and Lemma \ref{lm:free_upper_scaling}, we have
	\begin{align*}
		\bbE_x[\tau_{B_n}]
		& \ge \sum_{y \in B(x, bn)} {G_{B_n} (x, y)} \ge \sum_{y \in B(x, bn) \setminus \{x\}} {C_1 G(x, y)} \asymp \sum_{y \in B(x, bn) \setminus \{x\}} {g(\aps{x - y})} \\
		& \asymp \int_{1}^{bn} {g(r) r^{d - 1} \dif r} = \int_{1}^{bn} {\frac{1}{r \phi(r^{-2})} \dif r} = \frac{1}{\phi(n^{-2})} \int_{1}^{bn} {\frac{1}{r} \frac{\phi(n^{-2})}{\phi(r^{-2})} \dif r} \\
		& \ge \frac{1}{a_2 \phi(n^{-2}) n^{2 \gamma_2}} \int_{1}^{bn} {r^{2 \gamma_2 - 1} \dif r} = \frac{1}{2a_2 \gamma_2 \phi(n^{-2})} \UGL{b^{2 \gamma_2} - \frac{1}{n^{2 \gamma_2}}} \ge \frac{b^{2 \gamma_2}}{4 a_2 \gamma_2 \phi(n^{-2})},
	\end{align*}
	for $n$ large enough. Hence, we can conclude that $\bbE_x[\tau_{B_n}] \ge C_2 / \phi(n^{-2})$, for all $x \in B_{\frac{an}{2}}$, for $n$ large enough and for some $C_2 > 0$. As usual, we can adjust the constant to get the statement of this proposition for every $n \in \bbN$. Notice that this is true regardless of the dimension because here, we can always plug in $\gamma_2 = 1$.
\end{PF}

Now we want to find the upper bound for $\bbE_x [\tau_{B_n}]$.
\begin{LM}\label{lm:E[tau_Bn]<=}
	There exists constant $C_3 > 0$ such that for all $n \in \bbN$
	\begin{equation}\label{eq:E[tau_Bn]<=}
		\bbE_x[\tau_{B_n}] \le \frac{C_3}{\phi(n^{-2})}, \quad \forall\, x \in B_n.
	\end{equation}
\end{LM}
\begin{PF}
	We define the process $M^f = (M_n^f)_{n \ge 0}$ with
	\begin{equation*}
		M_n^f := f(X_n) - f(X_0) - \sum_{k = 0}^{n - 1} {(Af)(X_k)}
	\end{equation*}
	where $f$ is a function defined on $\bbZ^d$ with values in $\bbR$, $A$ is defined as in \eqref{eq:def_od_A} and $X = (X_n)_{n \ge 0}$ is a subordinate random walk. By \cite[Theorem 4.1.2]{No97}, the process $M^f$ is a martingale for every bounded function $f$. Let $f := \bbjedan_{B_{2n}}$ and $x \in B_n$. By the optional stopping theorem, we have
	\begin{equation*}
		\bbE_x [M_{\tau_{B_n}}^f] = \bbE_x \UGL{f(X_{\tau_{B_n}}) - f(X_0) - \sum_{k = 0}^{\tau_{B_n} - 1} {(Af)(X_k)}} = \bbE_x [M_0^f] = 0.
	\end{equation*}
	Hence
		\begin{equation}\label{eq:pom_rez_za_E[tau_Bn]}
		\bbE_x \UGL{f(X_{\tau_{B_n}}) - f(X_0)} = \bbE_x \UGL{\sum_{k = 0}^{\tau_{B_n} - 1} {(Af)(X_k)}}.
	\end{equation}
	We now investigate both sides of the relation \eqref{eq:pom_rez_za_E[tau_Bn]}. For every $k < \tau_{B_n}$, $X_k \in B_n$, and for every $y \in B_n$, using Proposition \ref{prop:asymptotics_for_jumps}, \eqref{eq:svojstvo_skaliranja} and \eqref{eq:upper_svojstvo_skaliranja}, we have
	\begin{align*}
		(Af)(y)
		& = \sum_{u \in \bbZ^d} {\bbP(X_1 = u - y) (f(u) - f(y))} \asymp -\sum_{u \in B_{2n}^c} {\aps{u - y}^{-d} \phi(\aps{u - y}^{-2})} \\
		& \asymp -\int_{n}^{\infty} {r^{-d} \phi(r^{-2}) r^{d - 1} \dif r} = -\phi(n^{-2}) \int_{n}^{\infty} {r^{-1} \frac{\phi(r^{-2})}{\phi(n^{-2})} \dif r} \\
		& \asymp -\phi(n^{-2}) \int_{n}^{\infty} {r^{-1} \frac{n^{2\gamma_i}}{r^{2\gamma_i}} \dif r} = -\phi(n^{-2}) n^{2\gamma_i} \frac{n^{-2\gamma_i}}{2\gamma_i} \asymp -\phi(n^{-2}).
	\end{align*}
	Using the above estimate, we get	\begin{equation}\label{eq:izraz_za_E[tau_Bn]}
		\bbE_x \UGL{\sum_{k = 0}^{\tau_{B_n} - 1} {(Af)(X_k)}} \asymp \bbE_x \UGL{-\sum_{k = 0}^{\tau_{B_n} - 1} {\phi(n^{-2})}} = -\phi(n^{-2}) \bbE_x[\tau_{B_n}].
	\end{equation}
	Using \eqref{eq:pom_rez_za_E[tau_Bn]}, \eqref{eq:izraz_za_E[tau_Bn]} and  $\bbE_x [f(X_{\tau_{B_n}}) - f(X_0)] = \bbP_x (X_{\tau_{B_n}} \in B_{2n}) - 1 = -\bbP_x(X_{\tau_{B_n}} \in B_{2n}^c)$, we get
	\begin{equation*}
		\bbP_x (X_{\tau_{B_n}} \in B_{2n}^c) \asymp \phi(n^{-2}) \bbE_x[\tau_{B_n}]
	\end{equation*}
	and this implies	\begin{equation}\label{eq:gornja_ograda_za_drugi_dio_prve_sume}
		\bbE_x [\tau_{B_n}] \le \frac{C_3 \bbP_x(X_{\tau_{B_n}} \in B_{2n}^c)}{\phi(n^{-2})} \le \frac{C_3}{\phi(n^{-2})}.
	\end{equation}
\end{PF}
In the next two results we develop estimates for the Green function of a ball. We define $A(r, s) := \{x \in \bbZ^d : r \le \aps{x} < s\}$ for $r, s \in \bbR$, $0 < r < s$.
\begin{PROP}\label{prop:upper_bound_for_G_Bn}
	There exists constant $C_4 > 0$ such that for all $n \in \bbN$	\begin{equation}\label{eq:upper_bound_for_G_Bn}
		G_{B_n}(x, y) \le C_4n^{-d} \eta(y), \quad \forall\, x \in B_{\frac{bn}{2}}, y \in A(bn, n),
	\end{equation}
	where $\eta(y) = \bbE_y[\tau_{B_n}]$, $b = a/2$, and $a \in \langle 0, 1/3 \rangle$ is as in Lemma \ref{lm:G_(B_n)>=G}.
\end{PROP}
\begin{PF}
	Let $x \in B_{\frac{bn}{2}}$ and $y \in A(bn, n)$. We define function $h(z) := G_{B_n} (x, z)$. Notice that for $z \in B_n \setminus \{x\}$ we have
	\begin{equation*}
		h(z) = G_{B_n}(x, z) = G_{B_n}(z, x) = \sum_{y \in \bbZ^d} {\bbP(X_1 = y - z) G_{B_n} (y, x)} = \sum_{y \in \bbZ^d} {\bbP(X_1 = y - z)} h(y).
	\end{equation*}
	Hence, $h$ is a harmonic function in $B_n \setminus \{x\}$. If we take $z \in B(x, bn/8)^c$ then $\aps{z - x} \ge bn/8 \ge 1$ for $n$ large enough. Using Lemma \ref{lm:g_is_almost_decreasing} and Theorem \ref{tm:asymptotics_of_G} we get
	\begin{equation*}
		g(bn/8) \ge a_2^{-1} g(\aps{z - x}) \asymp G(x, z) \ge G_{B_n}(x, z) = h(z).
	\end{equation*}
	Hence, $h(z) \le kg(bn/8)$ for $z \in B(x, bn/8)^c$ and for some constant $k > 0$. Notice that $A(bn, n) \subseteq B(x, bn/8)^c$, hence $y \in B(x, bn/8)^c$. Using these facts together with Proposition \ref{prop:asymptotics_for_jumps}, we have
	\begin{align*}
		A
		& (h \wedge kg(bn/8))(y) = A(h \wedge kg(bn/8) - h)(y) \\
		& = \sum_{v \in \bbZ^d} {\bbP(X_1 = v - y)} (h(v) \wedge kg(bn/8) - h(v) - h(y) \wedge kg(bn/8) + h(y)) \\
		& \asymp \sum_{v \in B(x, bn/8)} {j(\aps{v - y}) (h(v) \wedge kg(bn/8) - h(v))} \ge -\sum_{v \in B(x, bn/8)} {j(\aps{v - y}) h(v)} \\
		& \ge -\sum_{v \in B(x, bn/8)} {a_1^{-1} j(bn/8) h(v)} = -a_1^{-1} j(bn/8) \sum_{v \in B(x, bn/8)} {G_{B_n}(x, v)} \ge -a_1^{-1} j(bn/8) \eta(x),
	\end{align*}
	where in the last line we used Lemma \ref{lm:j_is_almost_decreasing} together with $\aps{v - y} \ge bn/8 \ge 1$ for $v \in B(x, bn/8)$ and for $n$ large enough. Using \eqref{eq:j(ar)_le_noSC} we get $j(bn/8) \le (b/8)^{-d-2} j(n)$. Hence, using \eqref{eq:E[tau_Bn]<=}, we have
	\begin{equation*}
		A(h \wedge kg(bn/8))(y) \ge -c_1 n^{-d} \phi\OBL{n^{-2}} \eta(x) \ge -c_1 n^{-d} \phi\OBL{n^{-2}} C_3\OBL{\phi\OBL{n^{-2}}}^{-1} = -c_2 n^{-d}
	\end{equation*}
	for some $c_2 > 0$. On the other hand, using \eqref{eq:g(ar)_le_LSC} and Proposition \ref{prop:E[tau_Bn]>=} we get
	\begin{align*}
		g(bn/8)
		& \le a_1^{-1} (bn/8)^{-d + 2\gamma_1} g(n) \le (a_1 C_2)^{-1} (bn/8)^{-d + 2\gamma_1} n^{-d} \eta(z) = c_3n^{-d} \eta(z), \quad \forall z \in B_{bn}.
	\end{align*}
	Now we define $C_4 := (c_2 \vee kc_3) + 1$ and using
	\begin{equation*}
		h(z) \wedge kg(bn/8) \le kg(bn/8) \le kc_3n^{-d} \eta(z)
	\end{equation*}
	we get
	\begin{equation*}
		C_4n^{-d} \eta(z) - h(z) \wedge kg(bn/8) \ge (C_4 - kc_3)n^{-d} \eta(z) \ge 0,  \quad \forall\, z \in B_{bn}
	\end{equation*}
	So, if we define $u(\cdot) := C_4n^{-d}\eta(\cdot) - h(\cdot) \wedge kg(bn/8)$, we showed that $u$ is non-negative function on $B_{bn}$. It is obvious that it vanishes on $B_n^c$ and for $y \in A(bn, n)$ we have
	\begin{equation*}
		(Au)(y) = C_4n^{-d} (A\eta)(y) - A(h \wedge kg(bn/8))(y) \le -C_4n^{-d} + c_2n^{-d} < 0.
	\end{equation*}
	Since $u \ge 0$ on $B_{bn}$ and $u$ vanishes on $B_n^c$, if $\inf_{y \in \bbZ^d}u(y) < 0$ then there would exist $y_0 \in A(bn, n)$ such that $u(y_0) = \inf_{y \in \bbZ^d}u(y)$. But then, by Proposition \ref{tm:maximum_principle}, $(Au)(y_0) \ge 0$ which is in contradiction with $(Au)(y) < 0$ for $y \in A(bn, n)$. Hence,
	\begin{equation*}
		u(y) = C_4n^{-d}\eta(y) - h(y) \wedge kg(bn/8) \ge 0, \quad \forall\, y \in \bbZ^d
	\end{equation*}
	and then, because $h(y) \le kg(bn/8)$ for $y \in A(bn, n)$ we get
	\begin{equation*}
		 G_{B_n}(x, y) = h(y) \le C_4n^{-d} \eta(y), \quad \forall\, x \in B_{\frac{bn}{2}},\, y \in A(bn, n).
	\end{equation*}
\end{PF}
Now we will prove a proposition that will give us the lower bound for the Green function of a ball. We use the fact that $\aps{B_n \cap \bbZ^d} \ge c n^d$ for some constant $c > 0$, where $\aps{\cdot}$ denotes the cardinality of a set.
\begin{PROP}\label{prop:lower_bound_for_G_Bn}
	There exist $C_5 > 0$ and $b \le a/4$ such that for all $n \in \bbN$	\begin{equation}\label{eq:lower_bound_for_G_Bn}
		G_{B_n}(x, y) \ge C_5n^{-d}\eta(y), \quad \forall\, x \in B_{bn}, y \in A(an/2, n),
	\end{equation}
	where $a$ is as in Lemma \ref{lm:G_(B_n)>=G} and $\eta(y) = \bbE_y[\tau_{B_n}]$. 
\end{PROP}
\begin{PF}
	Let $a \in \langle 0, 1/3 \rangle$ as in Lemma \ref{lm:G_(B_n)>=G}. Then there exists $C_1 > 0$
	\begin{equation}\label{eq:prop_lower_bound_for_G_Bn}
		G_{B_n} (x, v) \ge C_1G(x, v), \quad x, v \in B_{an}.
	\end{equation}
	From Proposition \ref{prop:upper_bound_for_G_Bn} it follows that there exists constant $C_4 > 0$ such that
	\begin{equation}\label{eq:prop_upper_bound_for_G_Bn}
			G_{B_n}(x, v) \le C_4n^{-d}\eta(v), \quad x \in B_{an/4}, v \in A(an/2, n).
	\end{equation}
	From Lemma \ref{lm:E[tau_Bn]<=} we have
	\begin{equation}\label{eq:prop_upper_bound_for_eta}
		\eta(v) \le \frac{C_3}{\phi\OBL{n^{-2}}}, \quad v\in B_n,
	\end{equation}
	for some constant $C_3 > 0$. By Theorem \ref{tm:asymptotics_of_G} and \eqref{eq:def_of_g} there exists $c_1 > 0$ such that $G(x) \ge c_1 g(\aps{x})$, $x \neq 0$. Now we take
	\begin{equation*}
		b \le \min\VIT{\frac{a}{4}, \OBL{\frac{C_1 c_1}{2 a_2^2 C_3 C_4}}^{\frac{1}{d - 2\gamma_2}}}
	\end{equation*}
	and fix it. Notice that $(C_1 c_1) / (a_2^2 C_3 b^{d - 2\gamma_2}) \ge 2C_4$. Let $x \in B_{bn}$, $v \in B(x, bn)$. Since $b \le a/4$ we have $x, v \in B_{an}$. We want to prove that $G_{B_n}(x, v) \ge 2C_4 n^{-d} \eta(v)$. We will first prove that assertion for $x \neq v$. In that case we have $1 \le \aps{x - v}$. Since $v \in B(x, bn)$, we have $\aps{x - v} \le bn$ so we can use \eqref{eq:prop_lower_bound_for_G_Bn}, Lemma \ref{lm:g_is_almost_decreasing} and \eqref{eq:g(ar)_ge_USC} to get
	\begin{align}\label{al:expression_for_G_Bn(x,v)}
		G_{B_n}(x, v) \ge C_1 G(x, v) \ge \frac{C_1 c_1}{a_2}g(bn) \ge \frac{C_1 c_1}{a_2^2 b^{d - 2\gamma_2}} g(n) \ge \frac{2C_3 C_4}{n^d \phi(n^{-2})}.
	\end{align}
	Using \eqref{eq:prop_upper_bound_for_eta} and \eqref{al:expression_for_G_Bn(x,v)}, we get $G_{B_n}(x, v) \ge 2C_4 n^{-d} \eta(v)$ for $x \neq v$. Now we will prove that $G_{B_n}(x, x) \ge 2C_4 n^{-d} \eta(x)$, for $x \in B_{bn}$ and for $n$ large enough. First note that
\begin{equation*}
		\lim_{n \to \infty} {n^d \phi(n^{-2})} = \lim_{n \to \infty} {n^d\ \frac{\phi(n^{-2})}{\phi(1)}} \ge \lim_{n \to \infty} {n^d \frac{1}{a_2 n^{2 \gamma_2}}} = \lim_{n \to \infty} {\frac{1}{a_2} n^{d - 2\gamma_2}} = \infty,
	\end{equation*}
	since $d - 2\gamma_2 > 0$. Therefore 
\begin{equation*}
	2 C_4 n^{-d} \eta(x) \le \frac{2 C_4 C_3}{n^d \phi(n^{-2})} \le 1 \le G_{B_n}(x, x)
\end{equation*}
for $n$ large enough. Hence,
	\begin{equation}\label{eq:c2n^-deta<=0.5G_Bn}
		C_4n^{-d} \eta(v) \le \frac{1}{2}G_{B_n}(x, v), \quad \forall\, x \in B_{bn}, v \in B(x, bn).
	\end{equation}
	Now we fix $x \in B_{bn}$ and define the function
	\begin{equation*}
		h(v) := G_{B_n}(x, v) \wedge \OBL{C_4n^{-d}\eta(v)}.
	\end{equation*}
	From \eqref{eq:c2n^-deta<=0.5G_Bn} we have $h(v) \le \frac{1}{2}G_{B_n}(x, v)$ for $v \in B(x, bn)$. Recall that $G_{B_n}(x, \cdot)$ is harmonic in $A(an/2, n)$. Using \eqref{eq:prop_upper_bound_for_G_Bn} we get $h(y) = G_{B_n}(x, y)$ for $y \in A(an/2, n)$. Hence, for $y \in A(an/2, n)$
	\begin{align}\label{al:(Ah)(y)-pocetak}
		(Ah)(y)
		& = A(h(\cdot) - G_{B_n}(x, \cdot))(y) \asymp \sum_{v \in \bbZ^d} {j(\aps{v - y}) \OBL{h(v) - G_{B_n}(x, v) - h(y) + G_{B_n}(x, y)}} \nonumber \\
		& \le \sum_{v \in B(x, bn)} {j(\aps{v - y})\OBL{h(v) - G_{B_n}(x, v)}} \le -(a_1/2) j(2n) \sum_{v \in B(x, bn)} {G_{B_n}(x, v)},
	\end{align}
	where we used Proposition \ref{prop:asymptotics_for_jumps} and Lemma \ref{lm:j_is_almost_decreasing} together with $1 \le \aps{v - y} \le 2n$. Using \eqref{al:expression_for_G_Bn(x,v)} and $\aps{B_n \cap \bbZ^d} \ge c_2 n^d$, we get
	\begin{equation}\label{eq:(Ah)(y)-pom1}
		\sum_{v \in B(x, bn)} {G_{B_n}(x, v)} \ge \frac{2C_3 C_4}{n^d \phi\OBL{n^{-2}}} \aps{B_{bn} \cap \bbZ^d} \ge \frac{2c_2 C_3 C_4}{n^d \phi\OBL{n^{-2}}} (bn)^d = \frac{c_3}{\phi\OBL{n^{-2}}}.
	\end{equation}
	Using \eqref{eq:j(ar)_ge_noSC} we get $j(2n) \ge 2^{-d - 2}j(n)$. When we put this together with \eqref{al:(Ah)(y)-pocetak} and \eqref{eq:(Ah)(y)-pom1}, we get
	\begin{equation*}
		(Ah)(y) \le -c_4n^{-d}.
	\end{equation*}
	Define $u := h - \kappa\eta$, where
	\begin{equation*}
		\kappa := \min\VIT{\frac{c_4}{2}, \frac{c_5}{2}, \frac{C_4}{2}}n^{-d},
	\end{equation*}
	where $c_5 > 0$ will be defined later. For $y \in A(an/2, n)$
	\begin{equation*}
		(Au)(y) = (Ah)(y) - \kappa(A\eta)(y) \le -c_4n^{-d} + \kappa \le -c_4n^{-d} + \frac{c_4}{2} n^{-d} = -\frac{c_4}{2}n^{-d} < 0.
	\end{equation*}
	For $x \in B_{bn} \subseteq B_{an/2}$, $v \in B_{an/2}$ we have $\aps{x - v} \le an \le n$. We will first assume that $x \neq v$ so that we can use Theorem \ref{tm:asymptotics_of_G}, Lemma \ref{lm:g_is_almost_decreasing} and \eqref{eq:g(ar)_ge_USC}. In this case, we have
	\begin{align*}
		G_{B_n}(x, v)
		& \ge C_1G(x, v) \asymp g(\aps{x - v}) \ge \frac{1}{a_2}g(an) \ge \frac{1}{a_2^2 a^{d - 2\gamma_2}} g(n) \ge \frac{1}{a_2^2 C_3 a^{d - 2\gamma_2}}n^{-d} \eta(v).
	\end{align*}
	So, $G_{B_n}(x, v) \ge c_5 n^{-d} \eta(v)$ for some constant $c_5 > 0$ and for $x \neq v$. If $x = v$ we can use the same arguments that we used when we were proving that $G_{B_n}(x, x) \ge 2C_4 n^{-d} \eta(x)$ for $n$ large enough to prove that $G_{B_n}(x, x) \ge c_5 n^{-d} \eta(x)$ for $n$ large enough. Hence, $G_{B_n}(x, v) \ge c_5 n^{-d} \eta(v)$ for all $x \in B_{bn}$ and $v \in B_{an/2}$ and for $n$ large enough. Now we have
	\begin{equation*}
		h(v) = G_{B_n}(x, v) \wedge \OBL{C_4n^{-d}\eta(v)} \ge \OBL{c_5n^{-d}\eta(v)} \wedge \OBL{C_4n^{-d}\eta(v)} = (C_4 \wedge c_5) n^{-d}\eta(v).
	\end{equation*}
	Hence,
	\begin{equation*}
		u(v) = h(v) - \kappa \eta(v) \ge (C_4 \wedge c_5) n^{-d} \eta(v) - \OBL{\frac{C_4}{2} \wedge \frac{c_5}{2}} n^{-d} \eta(v) \ge 0.
	\end{equation*}
	Since $u(v) \ge 0$ for $v \in B_{an/2}$, $u(v) = 0$ for $v \in B_n^c$ and $(Au)(v) < 0$ for $v \in A(an/2, n)$ we can use the same argument as in Proposition \ref{prop:upper_bound_for_G_Bn} to conclude by Proposition \ref{tm:maximum_principle} that $u(y) \ge 0$ for all $y \in \bbZ^d$. Since $G_{B_n}(x, y) \le C_4 n^{-d} \eta(y)$ for $x \in B_{an/4}, y \in A(an/2, n)$ we have $h(y) = G_{B_n}(x, y)$. Using that, we have
	\begin{equation*}
		G_{B_n}(x, y) \ge \kappa \eta(y) = C_5n^{-d}\eta(y), \quad x\in B_{bn}, y \in A(an/2, n),
	\end{equation*}
	for $n$ large enough. As before, we can change the constant and get \eqref{eq:lower_bound_for_G_Bn} for all $n \in \bbN$.
\end{PF}
Using last two propositions, we have the next corollary.
\begin{COR}\label{cor:G_Bn-upper_and_lower}
	Assume \eqref{eq:svojstvo_skaliranja} and \eqref{eq:upper_svojstvo_skaliranja}. Then there exist constants $C_6, C_7 > 0$ and $b_1, b_2 \in \langle 0, \frac{1}{2} \rangle$, $2b_1 \le b_2$ such that
	\begin{equation}\label{eq:G_Bn-upper_and_lower}
		C_6 n^{-d} \bbE_y[\tau_{B_n}] \le G_{B_n}(x, y) \le C_7 n^{-d} \bbE_y[\tau_{B_n}], \quad \forall\, x \in B_{b_1n}, y \in A(b_2n, n).
	\end{equation}
\end{COR}
\begin{PF}
	This corollary follows directly from Proposition \ref{prop:upper_bound_for_G_Bn} and Proposition \ref{prop:lower_bound_for_G_Bn}. We can set $b_2 = a/2$ where $a \in \langle 0, 1/3 \rangle$ is as in Lemma \ref{lm:G_(B_n)>=G} and $b_1 = b$ where $b \le a/4$ is as in Proposition \ref{prop:lower_bound_for_G_Bn}.
\end{PF}

\section{Proof of the Harnack inequality}\label{sec:Proof_of_Harnack_inequality}
\noindent
We start this section with the proof of the proposition that will be crucial for the remaining part of our paper.
\begin{PROP}\label{prop:Ikeda-Watanabe}
	Let $\FJADEF{f}{\bbZ^d \times \bbZ^d}{[0, \infty\rangle}$ be a function and $B \subset \bbZ^d$ a finite set. For every $x \in B$ we have
	\begin{equation}\label{eq:Ikeda-Watanabe}
		\bbE_x \UGL{f(X_{\tau_B - 1}, X_{\tau_B})} = \sum_{y \in B} {G_B (x, y) \bbE \UGL{f(y, y + X_1) \bbjedan_{\{y + X_1 \notin B\}}}}.
	\end{equation}
\end{PROP}
\begin{PF}
	We have
	\begin{equation*}
		\bbE_x \UGL{f(X_{\tau_B - 1}, X_{\tau_B})} = \sum_{y \in B, z \in B^c} {\bbP_x (X_{\tau_B - 1} = y, X_{\tau_B} = z) f(y, z)}.
	\end{equation*}
	Using \eqref{eq:X_as_RW}, we get
	\begin{align*}
		\bbP_x (X_{\tau_B - 1}
		& = y, X_{\tau_B} = z) = \sum_{m = 1}^{\infty} {\bbP_x (X_{\tau_B - 1} = y, X_{\tau_B} = z, \tau_B = m)} \\
		& = \sum_{m = 1}^{\infty} {\bbP (x + X_{m - 1} + \xi_m = z, x + X_{m - 1} = y, X_1, \ldots, X_{m - 2} \in B - x)} \\
		& = \sum_{m = 1}^{\infty} {\bbP (\xi_m = z - y) \bbP (x + X_{m - 1} = y, X_1, \ldots, X_{m - 2} \in B - x)} \\
		& = \bbP (\xi_1 = z - y) \sum_{m = 1}^{\infty} {\bbP_x (X_{m - 1} = y, X_1, \ldots, X_{m - 2} \in B)} \\
		& = \bbP (X_1 = z - y) \sum_{m = 1}^{\infty} {\bbP_x (X_{m - 1} = y, \tau_B > m - 1)} = \bbP (X_1 = z - y) G_B (x, y).
	\end{align*}
	Hence,
	\begin{align*}
		\bbE_x \UGL{f(X_{\tau_B - 1}, X_{\tau_B})}
		& = \sum_{y \in B, z \in B^c} {f(y, z) G_B (x, y) \bbP (y + X_1 = z)} \\
		& = \sum_{y \in B} {G_B (x, y) \bbE \UGL{f(y, y + X_1) \bbjedan_{\{y + X_1 \notin B\}}}}.
	\end{align*}
\end{PF}
\begin{REM}
	Formula \eqref{eq:Ikeda-Watanabe} can be considered as a discrete counterpart of the continuous-time Ikeda-Watanabe formula. We will refer to it as discrete Ikeda-Watanabe 
formula.
\end{REM}
It can be proved that if $\FJADEF{f}{\bbZ^d}{[0, \infty \rangle}$ is harmonic in $B$, with respect to $X$, then $\{f(X_{n \wedge \tau_B}) : n \ge 0\}$ is a martingale with respect to the natural filtration of $X$ (proof is the same as \cite[Proposition 1.4.1]{La96}, except that we have a non-negative instead of a bounded function). Using this fact, we can prove the following lemma.
\begin{LM}\label{lm:optional_stopping}
	Let $B$ be a finite subset of $\,\bbZ^d$. Then $\FJADEF{f}{\bbZ^d}{[0, \infty \rangle}$ is harmonic in $B$, with respect to $X$, if and only if $f(x) = \bbE_x[f(X_{\tau_B})]$ for every $x \in B$.
\end{LM}
\begin{PF}
	Let us first assume that $\FJADEF{f}{\bbZ^d}{[0, \infty\rangle}$ is harmonic in $B$, with respect to $X$. We take arbitrary $x \in B$. By the martingale property $f(x) = \bbE_x[f(X_{n \wedge \tau_B})]$, for all $n \ge 1$. First, by Fatou's lemma we have $\bbE_x[f(X_{\tau_B})] \le f(x)$ so $f(X_{\tau_B})$ is a $\bbP_x$-integrable random variable. Since $B$ is a finite set, we have $f \le M$ on $B$, for some constant $M > 0$, and $\bbP_x(\tau_B < \infty) = 1$. Using these two facts, we get
	\begin{equation*}
		f(X_{n \wedge \tau_B}) = f(X_n) \bbjedan_{\{n < \tau_B\}} + f(X_{\tau_B}) \bbjedan_{\{\tau_B \le n\}} \le M + f(X_{\tau_B}).
	\end{equation*}
	Since the right hand side is $\bbP_x$-integrable, we can use the dominated convergence theorem and we get
	\begin{equation*}
		f(x) = \lim_{n \to \infty} {\bbE_x[f(X_{n \wedge \tau_B})]} = \bbE_x[\lim_{n \to \infty} {f(X_{n \wedge \tau_B})}] = \bbE_x[f(X_{\tau_B})].
	\end{equation*}
	
	On the other hand, if $f(x) = \bbE_x[f(X_{\tau_B})]$, for every $x \in B$, then for $x \in B$ we have
	\begin{equation*}
		f(x) = \sum_{y \in \bbZ^d} {\bbE_x\,[\,f(X_{\tau_B}) \mid X_1 = y\,]\, \bbP_x(X_1 = y)} = \sum_{y \in \bbZ^d} {p(x, y) \bbE_y[f(X_{\tau_B})]} = \sum_{y \in \bbZ^d} {p(x, y) f(y)}.
	\end{equation*}
\end{PF}
Hence, if we take $B \subset \bbZ^d$ finite and $\FJADEF{f}{\bbZ^d}{[0, \infty \rangle}$ harmonic in $B$, with respect to $X$, then by Lemma \ref{lm:optional_stopping} and the discrete Ikeda-Watanabe formula, we get
\begin{equation}\label{eq:pocetak_raspisa_h_pomocu_poissonove_jezgre}
	f(x) = \bbE_x \UGL{f(X_{\tau_B})} = \sum_{y \in B} {G_B (x, y) \bbE \UGL{f(y + X_1) \bbjedan_{\{y + X_1 \notin B\}}}}.
\end{equation}
Let us define the discrete Poisson kernel of a finite set $B \subset \bbZ^d$ by
\begin{equation}\label{eq:def_of_Poisson_kernel}
	K_B(x, z):= \sum_{y \in B} {G_{B} (x, y) \bbP(X_1 = z - y)} , \quad x \in B, z \in B^c.
\end{equation}
If the function $f$ is non-negative and harmonic in $B_n$, with respect to $X$, from \eqref{eq:pocetak_raspisa_h_pomocu_poissonove_jezgre} we have
\begin{align}\label{al:zavrsetak_raspisa_h_pomocu_poissonove_jezgre}
	f(x)
	& = \sum_{y \in B_n} {G_{B_n} (x, y) \sum_{z \notin B_n} {\bbE \UGL{f(y + X_1) \bbjedan_{\{y + X_1 \notin B_n\}} \mid X_1 = z - y} \bbP(X_1 = z - y)}} \nonumber \\
	& = \sum_{z \notin B_n} {\sum_{y \in B_n} {G_{B_n} (x, y) \bbE \UGL{f(y + z - y) \bbjedan_{\{y + z - y \notin B_n\}}} \bbP(X_1 = z - y)}} \nonumber \\
	& = \sum_{z \notin B_n} {f(z) \OBL{\sum_{y \in B_n} {G_{B_n} (x, y) \bbP(X_1 = z - y)}}} = \sum_{z \notin B_n} {f(z) K_{B_n}(x, z)}.
\end{align}
Now we are ready to show that the Poisson kernel $K_{B_n}(x, z)$ is comparable to an expression that is independent of $x$. When we prove that, Harnack inequality will follow immediately.
\begin{LM}\label{lm:K_Bn_comparable_to_l}
	Assume \eqref{eq:svojstvo_skaliranja} and let $b_1, b_2 \in \langle 0, \frac{1}{2} \rangle$ be as in Corollary \ref{cor:G_Bn-upper_and_lower}. Then $K_{B_n}(x, z) \asymp l(z)$ for all $x \in B_{b_1n}$, where
	\begin{equation*}
		l(z) =  \frac{j(\aps{z})}{\phi\OBL{n^{-2}}} + n^{-d} \sum_{y \in A(b_2n, n)}{\bbE_y[\tau_{B_n}] j(\aps{z - y})}.
	\end{equation*}
\end{LM}
\begin{PF}
	Splitting the expression \eqref{eq:def_of_Poisson_kernel} for the Poisson kernel in two parts and using Proposition \ref{prop:asymptotics_for_jumps}, we get
	\begin{equation*}
		K_{B_n}(x, z) \asymp \sum_{y \in B_{b_2n}} {G_{B_n}(x, y) j(\aps{z - y})} + \sum_{y \in A(b_2n, n)} {G_{B_n}(x, y) j(\aps{z - y})}.
	\end{equation*}
	Since $G_{B_n}(x, y) \asymp n^{-d}\bbE_y[\tau_{B_n}]$ for $x \in B_{b_1n}$, $y \in A(b_2n, n)$, for the second sum in the upper expression we have
	\begin{equation}\label{eq:usporedivost_druge_sume}
		\sum_{y \in A(b_2n, n)} {G_{B_n}(x, y) j(\aps{z - y})} \asymp n^{-d} \sum_{y \in A(b_2n, n)}{\bbE_y[\tau_{B_n}] j(\aps{z - y})}.
	\end{equation}
	Now we look closely at the expression $\sum_{y \in B_{b_2n}} {G_{B_n}(x, y) j(\aps{z - y})}$. Using the fact that $y \in B_{b_2n}$, $b_2 \in \langle 0, \frac{1}{2} \rangle$ and $\aps{z} \ge n$ because $z \in B_n^c$, we have
	\begin{equation*}
		\aps{z - y} \le \aps{z} + \aps{y} \le \aps{z} + b_2n \le \aps{z} + b_2\aps{z} \le (1 + b_2)\aps{z} \le 2\aps{z}.
	\end{equation*}
	On the other hand
	\begin{equation}\label{eq:|z-y|<=z}
		\aps{z} \le \aps{z - y} + \aps{y} \le \aps{z - y} + b_2n \le \aps{z - y} + b_2\aps{z}.
	\end{equation}
	Hence,
	\begin{equation}\label{eq:|z-y|>=z}
		\frac{1}{2}\aps{z} \le (1 - b_2)\aps{z} \le \aps{z - y}.
	\end{equation}
	Combining \eqref{eq:|z-y|<=z}, \eqref{eq:|z-y|>=z} and using Lemma \ref{lm:j_is_almost_decreasing}, we have
	\begin{equation*}
		\frac{1}{a_1} j\OBL{\frac{1}{2}\aps{z}} \ge j(\aps{z - y}) \ge a_1 j(2\aps{z}).
	\end{equation*}
	Using \eqref{eq:j(ar)_le_noSC}, we get $j(\frac{1}{2} \aps{z}) \le 2^{d + 2}j(\aps{z}) = c_1 j(\aps{z})$. Similarly, from \eqref{eq:j(ar)_ge_noSC}, we get $j(2\aps{z}) \ge 2^{-d -2} j(\aps{z}) = c_2 j(\aps{z})$. Hence, $a_1 c_2 j(\aps{z}) \le a_1 j(2\aps{z}) \le j(\aps{z - y}) \le a_2^{-1} j\OBL{\frac{1}{2}\aps{z}} \le a_2^{-1} c_1 j(\aps{z})$ for some $c_1, c_2 > 0$. Therefore,
	\begin{equation}\label{eq:j(|z-y|)=j(|z|)}
		j(\aps{z - y}) \asymp j(\aps{z}), \quad y \in B_{b_2n},\, z \in B_n^c.
	\end{equation}
	Using \eqref{eq:j(|z-y|)=j(|z|)} we have
	\begin{equation*}
		\sum_{y \in B_{b_2n}} {G_{B_n}(x, y) j(\aps{z - y})} \asymp \sum_{y \in B_{b_2n}} {G_{B_n}(x, y) j(\aps{z})} = j(\aps{z}) \sum_{y \in B_{b_2n}} {G_{B_n}(x, y)}.
	\end{equation*}
	Now we want to show that $\sum_{y \in B_{b_2n}} {G_{B_n}(x, y)} \asymp 1 / \phi\OBL{n^{-2}}$. Using the fact that $G_{B_n}$ is non-negative function and that $\bbE_x[\tau_{B_n}] \le C_3 / \phi\OBL{n^{-2}}$ for $x \in B_n$ we have
	\begin{equation}\label{eq:sum_G_Bn<=1/phi}
		\sum_{y \in B_{b_2n}} {G_{B_n}(x, y)} \le \sum_{y \in B_n} {G_{B_n}(x, y)} = \bbE_x[\tau_{B_n}] \le \frac{C_3}{\phi\OBL{n^{-2}}}.
	\end{equation}
	To prove the other inequality we will use Lemma \ref{lm:G_(B_n)>=G}, Theorem \ref{tm:asymptotics_of_G}, Lemma \ref{lm:g_is_almost_decreasing} together with $1 \le \aps{x - y} \le 2b_2n$, $\aps{B_n \cap \bbZ^d} \ge c_3n^d$ and Lemma \ref{lm:free_upper_scaling}. Thus
	\begin{align*}
		\sum_{y \in B_{b_2n}} {G_{B_n}(x, y)}
		& \ge C_1\sum_{y \in B_{b_2n} \setminus \{x\}} {G(x, y)} \asymp \sum_{y \in B_{b_2n} \setminus \{x\}} {g(\aps{x - y})} \\
		& \ge \frac{1}{a_2} (\aps{B_{b_2n} \cap \bbZ^d} - 1) g(2b_2n) \ge \frac{1}{a_2} \frac{c_3}{2}(b_2n)^d \frac{1}{2^d (b_2n)^d} \frac{1}{\phi\OBL{n^{-2}}} \frac{\phi(n^{-2})}{\phi((2b_2n)^{-2})} \\
		& \ge \frac{c_3}{2 a_2} \frac{1}{2^d \phi(n^{-2})} (2b_2)^2 \ge \frac{c_3 (2b_2)^2}{2^{d + 1} a_2} \frac{1}{\phi\OBL{n^{-2}}}.
	\end{align*}
	Hence,
	\begin{equation}\label{eq:sum_G_Bn>=1/phi}
		\sum_{y \in B_{b_2n}} {G_{B_n}(x, y)} \ge \frac{c_4}{\phi\OBL{n^{-2}}}.
	\end{equation}
	From \eqref{eq:sum_G_Bn<=1/phi} and \eqref{eq:sum_G_Bn>=1/phi} we have
	\begin{equation}\label{eq:sum_G_Bn=1/phi}
		\sum_{y \in B_{b_2n}} {G_{B_n}(x, y)} \asymp \frac{1}{\phi\OBL{n^{-2}}}.
	\end{equation}
	Finally, using \eqref{eq:j(|z-y|)=j(|z|)} and \eqref{eq:sum_G_Bn=1/phi} we have
	\begin{equation}\label{eq:usporedivost_prve_sume}
		\sum_{y \in B_{b_2n}} {G_{B_n}(x, y) j(\aps{z - y})} \asymp \frac{j(\aps{z})}{\phi\OBL{n^{-2}}}.
	\end{equation}
	And now, from \eqref{eq:usporedivost_prve_sume} and \eqref{eq:usporedivost_druge_sume} we have the statement of the lemma.
\end{PF}
Lemma \ref{lm:K_Bn_comparable_to_l} basically states that there exist constants $C_8, C_9 > 0$ such that
\begin{equation}\label{eq:K_Bn_omedjeno_s_l(z)}
	C_8 l(z) \le K_{B_n}(x, z) \le C_9 l(z), \quad x \in B_{b_1n}, \, z \in B_n^c.
\end{equation}
Now we are ready to prove our main result.
\begin{PF}[Proof of Theorem \ref{tm:Harnack_inequality}]
	Notice that, because of the spatial homogeneity, it is enough to prove this result for balls centered at the origin. We will prove the theorem for $a = b_1$, where $b_1$ is as in Corollary \ref{cor:G_Bn-upper_and_lower}. General case follows using the standard Harnack chain argument. Let $x_1, x_2 \in B_{b_1n}$. Using \eqref{eq:K_Bn_omedjeno_s_l(z)} we get
	\begin{equation*}
		K_{B_n}(x_1, z) \le C_9 l(z) = \frac{C_9}{C_8} C_8 l(z) \le \frac{C_9}{C_8} K_{B_n}(x_2, z).
	\end{equation*}
	Now we can multiply both sides with $f(z) \ge 0$ and sum over all $z \notin B_n$ and we get
	\begin{equation*}
		\sum_{z \notin B_n} {f(z) K_{B_n}(x_1, z)} \le \frac{C_9}{C_8} \sum_{z \notin B_n} {f(z) K_{B_n}(x_2, z)}.
	\end{equation*}
	If we look at the expression \eqref{al:zavrsetak_raspisa_h_pomocu_poissonove_jezgre} we see that this means
	\begin{equation*}
		f(x_1) \le \frac{C_9}{C_8} f(x_2)
	\end{equation*}
	and that is what we wanted to prove.
\end{PF}

\noindent
\textbf{Acknowledgement:} This work has been supported in part by Croatian Science Foundation under the project 3526.

\bibliographystyle{babamspl}
\bibliography{Harnack}

\vspace{20px}

Ante Mimica, 20-Jan-1981 - 9-Jun-2016, \url{https://web.math.pmf.unizg.hr/~amimica/}

\vspace{10px}
Stjepan Šebek, Faculty of Electrical Engineering and Computing, University of Zagreb, 10000 Zagreb, Croatia

\textit{E-mail address}: stjepan.sebek@fer.hr
\end{document}